\documentclass[12pt,a4paper]{article}

\usepackage{amsmath}    % need for subequations
\usepackage{amssymb}
\usepackage{graphicx}   % need for figures
\usepackage{verbatim}   % useful for program listings
\usepackage{color}      % use if color is used in text
\usepackage{hyperref}   % use for hypertext links, including those to external documents and URLs
\usepackage{enumerate}
\usepackage{listings}
\usepackage{xcolor}
\usepackage{slashed}
\usepackage{amsthm}
\usepackage{amsmath}    % need for subequations
\usepackage{amssymb}
\usepackage{mathrsfs}
\usepackage{bbold}
\usepackage{lscape}
\usepackage{caption}

\usepackage[final]{showlabels}

\usepackage{amsthm}
\newtheorem{definition}{Definition}[section]

\newtheorem{theorem}{Theorem}[section]

\newtheorem{notation}{Notation}[section]

\newcommand{\del}{\partial}
\renewcommand{\theta}{\vartheta}
\renewcommand{\phi}{\varphi}

\newcommand{\dd}{\mathrm{d}}

\newcommand{\ii}{\mathbb{i}}

\renewcommand{\and}{\wedge}

\newcommand{\sinc}{\mathrm{sinc\,}}

\newcommand{\sign}{\mathrm{sign}\,}
\newcommand{\diag}{\mathrm{diag}\,}

\newcommand{\order}{N}

\renewcommand{\title}{Extensions of Active Flux to arbitrary order of accuracy}

\newcommand{\authorOne}{Rémi Abgrall\footnote{Institute for Mathematics and Computational Science, University of Zurich, Winterthurerstrasse 190, 8057 Zurich, Switzerland}}
\newcommand{\authorTwo}{Wasilij Barsukow\footnote{Institute for Mathematics, CNRS UMR 5251, University of Bordeaux, 351 Cours de la Libération, 33405 Talence, France}}
\begin{document}

\begin{center} \Large
\title

\vspace{1cm}

\date{}
\normalsize

\authorOne, \authorTwo
\end{center}

\begin{abstract}

% #############################################################################################################################
Active Flux is a recently developed numerical method for hyperbolic conservation laws. Its classical degrees of freedom are cell averages and point values at cell interfaces. These latter are shared between adjacent cells, leading to a globally continuous reconstruction. The update of the point values includes upwinding, but without solving a Riemann Problem. The update of the cell average requires a flux at the cell interface, which can be immediately obtained using the point values. This paper explores different extensions of Active Flux to arbitrarily high order of accuracy, while maintaining the idea of global continuity. We propose to either increase the stencil while keeping the same degrees of freedom, or to increase the number of point values, or to include higher moments as new degrees of freedom. These extensions have different properties, and reflect different views upon the relation of Active Flux to the families of Finite Volume, Finite Difference and Finite Element methods. 
% #############################################################################################################################



\end{abstract}

\section{Introduction}

Solutions to the initial value problem for hyperbolic $m\times m$ systems of conservation laws
\begin{align}
 \del_t q + \del_x f(q) &= 0 & q &\colon \mathbb R^+_0 \times \mathbb R \to \mathbb R^m \label{eq:conslaw}\\
 && f &\colon \mathbb R^m \to \mathbb R^m
\end{align}
generically develop discontinuities in finite time, even if the initial data are smooth. The natural setting therefore are weak solutions. For convergence to the weak solution of \eqref{eq:conslaw}, a numerical method needs to be conservative (Lax-Wendroff theorem).

A popular way to derive numerical methods (due to Godunov \cite{godunov59difference}) is to introduce discontinuities at every cell interface (reconstruction step), and to evolve piecewise constant data over a short period of time (Riemann solver). Even if a more complicated reconstruction is used inside the cell, or if an approximate evolution is used instead of the exact one, the main idea of Godunov's approach remains to introduce discontinuities at every cell interface. This approach has also been adapted to Finite Element methods, which has led to the development of Discontinuous Galerkin (DG) methods (e.g. \cite{lesaint74,cockburn91first,cockburn1989a,cockburn89}). These methods are conservative.

However, convergence to the weak solution is not enough in practice. In view of the big computational effort associated with grid refinement (particularly in multi-d), there is ongoing interest in guaranteeing properties of numerical solutions for coarse grids already. Some of them are essential for the simulation to continue running, e.g. non-negativity of the density. Other shortcomings might not cause a simulation to crash, but still are a big problem in practice, in particular if they require the computational grid to be much finer than reasonable. For example it might turn out to be necessary to resolve a vortex with a number of grid cells that depends on how long one intends to run the simulation for. This can occur in simulations of fluid flow phenomena that are smooth (low Mach number/incompressible limit). They are not well approximated by standard Godunov methods on coarse grids because the effects of numerical diffusion become more and more pronounced as the simulation runs. Grid refinement increases the time scales on which diffusion becomes overwhelming, but -- besides slowing down the simulation -- does not solve the problem in a fundamental way. 

The difficulties of the Godunov method in regions of smooth flow are not surprising, as it introduces shocks everywhere in the computational domain. For example, the failure of the Godunov method to resolve the low Mach number limit on coarse grids has been traced back directly to the acoustic waves introduced at every time step at every cell interface (in \cite{guillard04}). For low Mach number flow, modifications of Godunov methods have been suggested which restore the correct behaviour, but they are ad hoc and cannot guarantee adequate treatment of other aspects of smooth multi-dimensional flow (for example vortices). Note at this point that in the quest for a new numerical method we expect it to be able to work well in presence of shocks as well as when they are absent.

The argumentation of \cite{guillard04} assumes that the multi-dimensional numerical method is constructed using a one-dimensional Riemann solver. The fluxes through each interface are computed while ignoring the influence of the corners of the cell. Such an approach shall be called directionally-split. Truly multi-dimensional phenomena, where contributions from different directions balance each other, tend to be a challenge for dimensionally split methods. The low Mach number/incompressible limit is one example of such a truly multi-dimensional phenomenon, because the divergence-free condition is trivial in one spatial dimension. One might therefore conjecture that the difficulties do not originate in the Godunov method per se, but in the fact that quasi-one-dimensional solvers are applied in a directionally split way, and thus multi-dimensionality is not taken into account appropriately. One might think that perhaps it would suffice to take into account all the multi-dimensional interactions in multi-dimensional Riemann problems to obtain a Godunov-type method that resolves smooth multi-dimensional problems accurately on coarse grids.

Unfortunately, the situation is not that simple. In \cite{barsukow17} this line of thought has been developed completely for the linear acoustic equations, which have interesting multi-dimensional features that are not captured by standard dimensionally-split Godunov-type methods on coarse grids (e.g. an analogue of the low Mach number limit, and the involution of vorticity). The complete solution of the four-quadrant Riemann problem for linear acoustics has been obtained and a two-dimensional Godunov method on Cartesian grids constructed in \cite{barsukow17}, whose evolution step is exact and fully multi-dimensional. None of the desirable properties (low Mach number compliance, vorticity preservation) was found to hold true for this method. Taking into account all multi-dimensional interactions from multi-dimensional Riemann problems is not enough (and, of course, quite an effort).

Similar observations described in \cite{roe16,roe17a} have sparked the development of the Active Flux method by Roe and collaborators (\!\!\!\cite{eymann11}, based on \cite{vanleer77}), a numerical method with a continuous reconstruction, whose evolution step thus requires a short-time solution of the IVP for \eqref{eq:conslaw} with \emph{continuous} data. 

This does not mean that this method is only applicable in the subsonic regime. Continuous data can self-steepen over the course of a time step, and in such cases, the update needs to be performed carefully, before the discontinuity is projected onto a continuous reconstruction again. In a typical simulation, however, shocks are located on sets of codimension 1 (e.g. in 1-d on countably many points), and thus in this sense almost everywhere in the computational domain the assumption of continuity is correct. Compare this to the approach of a Godunov method solving a rarefaction: instead of maintaining the continuous solution, the rarefaction is cut apart in a staircase shape at every time step, and over the course of the evolution step every jump of the staircase is evolving again into a small rarefaction.

Leaving aside conceptual questions, from a practical point of view, the answers to the following three questions matter most: Is this method stable, possibly even under an explicit time integration? Does the numerical solution obtained with this method converge to the weak solution of \eqref{eq:conslaw}, i.e. is there an analogue of the Lax-Wendroff theorem? Does this method have favorable properties on smooth, multi-dimensional solutions? The answer to all three questions is yes, which makes it an interesting alternative to Godunov-type methods.

The degrees of freedom of the traditional Active Flux (\!\!\!\cite{vanleer77,eymann11}) are cell averages and, additionally, point values located at cell interfaces and shared by adjacent cells. This allows for a continuous (but not necessarily differentiable) reconstruction, which passes through the given point values, and whose average agrees with the given one. The evolution of the point value requires the solution of the IVP for \eqref{eq:conslaw} over a short time step for continuous initial data. This building block replaces the Riemann solver in Godunov-type methods. The evolution of the average is obtained as usual by integrating the conservation law over the cell and over a time step, and by applying Gauss' theorem. The fluxes at cell interfaces can be immediately evaluated using the available point values. The evolution of the averages therefore is conservative, and an analogue of a Lax-Wendroff theorem for Active Flux has been shown in \cite{abgrall20}. Active Flux decouples the problem of average evolution from the computation of interface values, which are declared independent degrees of freedom, hence the name Active Flux, although it would be more accurate to say ``Active Point Value at Cell Interface''. The additional point values allow for more flexibility. The update of the point values is not subject to the constraint of conservation, because conservation applies only to averages. Therefore it is possible to even switch variables and update the point values in, say, primitive variables. This has been pointed out in \cite{abgrall20}. A more thorough description to the traditional Active Flux method is given in Section \ref{sec:reviewtraditional}, as well as in \cite{eymann11,eymann13,barsukow18activeflux,kerkmann18,barsukow19activeflux}. 

The development of an arbitrary-order method requires chiseling out the distinctive features of the low order method. The aim of the present paper therefore is not only to propose practical high-order extensions of Active Flux. These extensions are tied to possible interpretations of Active Flux that allow to place it amidst other families of numerical methods. The three extensions that we propose allow to see Active Flux either as a coupled Finite Volume/Finite Difference method (Section \ref{sec:findiff}), or as an enriched Finite Volume method (Section \ref{sec:points}), or as a coupled Finite Element/Finite Difference method (Section \ref{sec:moments}). Which of the three interpretations, or maybe another, will turn out to be the most fruitful one, only future can tell. Therefore at this stage we content ourselves with developing them to a level that allows to solve the Euler equations, and with stating their respective advantages and disadvantages, whereby we cannot single out any of the three extensions as overperforming the others in all respects.

%In this paper we intend to investigate three different ways of extending Active Flux to arbitrary order of accuracy. We first note that high order of accuracy can be either reached by extending the stencil of finite difference operators, or by incorporating more degrees of freedom per cell. We follow both paths and achieve the former by extending a finite-difference version of Active Flux from \cite{abgrall20}. As for the latter, it is possible to increase the number of pointwise degrees of freedom, or to observe that an average is a 0-th moment, and to take into consideration higher moments to achieve higher order of accuracy. This paper thus presents and compares three different extensions of the traditional Active Flux to higher order of accuracy.

%The different state of development of the three suggestions we make is due to their different nature. This often is antagonized by the properties of the resulting. For example, the best CFL conditions are obtained for the most sophisticated evolution of the point values, which naturally is less well developed. We hope, however, that such desirable properties will spark research to fill the gaps.

This paper focuses on the one-dimensional case and defers the case of multiple dimensions to a forthcoming publication in order to keep the size of the paper reasonable. Additionally, in multiple spatial dimensions there are further aspects that need to be taken into account, and aspects of the numerical discretization are also different for Cartesian grids and unstructured ones. 

Here, the grid is assumed equidistant with cells $[x_{i-\frac12}, x_{i+\frac12}]$ where $x_{i\pm\frac12} = x_i \pm \frac{\Delta x}{2}$, $i\in \mathbb Z$. As customary, indices referring to time are denoted as superscripts, i.e. $t^n$ is the $n$-th time level. $C^k$ denotes the class of $k$ times continuously differentiable functions, and $P^k$ the set of polynomials of degree less or equal to $k$.

\section{Review: Third order Active Flux} \label{sec:reviewtraditional}

In this section, an overview of the existing Active Flux method(s) in one spatial dimension shall be given. Only third-order methods have been suggested so far in the literature. Nevertheless, to facilitate later extensions, we would like to start pointing out general features of Active Flux here already.

Active Flux has two distinctive features. First, the mixed type of degrees of freedom: a cell average, reminiscent of Finite Volume methods, and point values located at cell boundaries. The numerical solution is thus given by the two sets

\begin{align}
 &\{ \bar q_i^n \in \mathbb R^m \}_{i \in I \subset \mathbb Z} & &\{ q_{i+\frac12}^n \in \mathbb R^m \}_{i\in I \subset \mathbb Z}  \label{eq:afdofs}
\end{align}
with the interpretations
\begin{align}
 \bar q_i^n &\simeq \frac{1}{\Delta x} \int_{x_{i-\frac12}}^{x_{i+\frac12}} \dd x \, q(t^n, x) & q_{i+\frac12}^n &\simeq q(t^n, x_{i+\frac12})
\end{align}
The same formulae are used to initialize the degrees of freedom $\{ \bar q^0_i \}_i, \{ q_{i+\frac12}^0 \}_i$ given initial data $q_0 \colon \mathbb R \to \mathbb R^m$. Traditional degrees of freedom in Active Flux are one cell average per cell, and one (shared) point value located at every cell interface.

Second, the point values are shared between adjacent cells, contrary to the approach of e.g. DG methods. The point values at cell interfaces allow to immediately write down a conservative update of the cell average. Indeed, integrating the conservation law \eqref{eq:conslaw} over a computational cell $[x_{i-\frac12}, x_{i+\frac12}]$ and applying Gauss' theorem yields
\begin{align}
 \frac{\dd}{\dd t} \left( \frac{1}{\Delta x} \int_{x_{i-\frac12}}^{x_{i+\frac12}} \dd x \, q(t, x)  \right) + \frac{f(q(t, x_{i+\frac12})) - f(q(t, x_{i-\frac12}))}{\Delta x} &= 0 \label{eq:averageupdategeneral}
\end{align}

\subsection{Active Flux with an evolution operator} \label{ssec:afwithevolutionop}

The traditional, earliest version of Active Flux (\cite{vanleer77,eymann11}) makes use of an evolution operator, i.e. of an operator that maps initial data to the exact solution at time $t$. It was later relaxed to allow for an approximate time evolution.

\begin{definition}[Traditional Active Flux] \label{def:activefluxtraditional}
Inside cell $i$, consider a reconstruction 
\begin{align}
q_{\text{recon},i} &\in C^0 \cap L^1_\text{loc} & q_{\text{recon},i} &\colon (\mathbb R^m)^3 \times \left[-\frac{\Delta x}{2},\frac{\Delta x}{2}\right] \to \mathbb R^m \label{eq:recon3rdstandard}
\end{align}
with the properties
\begin{align}
q_{\text{recon},i}\left(q_{i-\frac12}, \bar q_i, q_{i+\frac12}, \pm \frac{\Delta x}{2}\right ) &= q_{i\pm\frac12}\\
 \frac{1}{\Delta x} \int_{-\frac{\Delta x}{2}}^{\frac{\Delta x}{2}} \dd x \, q_{\text{recon},i}\left(q_{i-\frac12}, \bar q_i, q_{i+\frac12}, x\right ) &= \bar q_i 
 %q_{\text{recon},i}\left(\phi(x_{i-\frac12}), \frac{1}{\Delta x} \int_{x_{i-\frac12}}^{x_{i+\frac12}} \dd x \,\phi(x), \phi(x_{i+\frac12}), x\right ) = \phi(x_i + x) + \mathcal O(\Delta x^N) \qquad \forall x \in \left[-\frac{\Delta x}{2},\frac{\Delta x}{2}\right]
\end{align}
%for any analytic function $\phi \colon \mathbb R \to \mathbb  R$ and $N \geq 1$. 
Define the global reconstruction 
\begin{align}
q_\text{recon} &\in C^0 \cap L^1_\text{loc} & q_\text{recon} &\colon \mathbb R \to \mathbb R^m
\end{align}
\begin{align}
 q_\text{recon}(x) &:= q_{\text{recon},i}\left(q_{i-\frac12}, \bar q_i, q_{i+\frac12}, q_{i,1}, \ldots, q_{i,k}, x-x_i\right ) \quad \text{ if }x\in [x_{i-\frac12},x_{i+\frac12}]
\end{align}

The \emph{traditional Active Flux method} is the following semi-discretization
\begin{align}
 \begin{cases} \displaystyle \frac{\dd}{\dd t} \bar q_i(t) = - \frac{f(q_{i+\frac12}(t)) - f(q_{i-\frac12}(t))}{\Delta x} \\ \\
 \displaystyle q_{i+\frac12}(t) = \Big ( \text{solution at $x=x_{i+\frac12}$ of the IVP \eqref{eq:conslaw} with initial data } q_{\text{recon}} \Big ) + \mathcal O(t^3) \end{cases} \label{eq:activefluxdeftraditional}
 %\frac{\dd}{\dd t} q_{i+\frac12}(t) &= -R\Big(q_{i+\frac12-k}(t), \bar q_{i-k+1}(t), q_{i+\frac12-k+1}(t), \ldots, \bar q_{i + m}(t), q_{i+\frac12+m}(t)\Big)  \qquad i \in \mathbb Z, k\geq 0, m \geq 0 \label{eq:activefluxdef2}
\end{align}
of \eqref{eq:conslaw} with the interpretations
\begin{align}
 \bar q_i(t) &\simeq \frac{1}{\Delta x} \int_{x_{i-\frac12}}^{x_{i+\frac12}} \dd x \, q(t, x)  &
 q_{i+\frac12}(t) &\simeq q(t, x_{i+\frac12}) \label{eq:traditionaldofinterpretations}
\end{align}

%and $R$ a consistent approximation of $\del_x f(q)$.
\end{definition}

In the following, we often drop the dependence of $q_{\text{recon},i}$ on its parameters other than $x$. When the time step is important, we denote by $q_{\text{recon},i}^n$ the reconstruction obtained by using values of time step $n$.

Usually, $q_{\text{recon},i} \in P^2$ but limiting might require other choices, discussed in Section \ref{ssec:limiting}. An Active Flux method with this choice of reconstruction is at most third order accurate. 

Observe that the availability of an evolution operator in \eqref{eq:activefluxdeftraditional}, i.e. the possibility to solve the IVP at least approximately, allows to immediately rewrite the update of the averages as an explicit one by integrating over the time step:
\begin{align}
\frac{\bar q_i^{n+1} - \bar q_i^n}{\Delta t} = - \frac{\frac{1}{\Delta t} \int_{t^n}^{t^{n+1}} \dd t \, f(q_{i+\frac12}(t)) - \frac{1}{\Delta t} \int_{t^n}^{t^{n+1}} \dd t \,  f(q_{i-\frac12}(t))}{\Delta x} 
\end{align}
This equation is exact, i.e. its order of accuracy is entirely given by the order of accuracy of the point values. In order to use it in practice, the integrals are replaced by quadrature. First, approximations
\begin{align}
 q_{i+\frac12}^{n+\frac{\ell}{2}} = q_{i+\frac12}\left(  \ell \frac{\Delta t}{2} \right) + \mathcal O(\Delta t^3) \qquad \ell = 1,2
\end{align}
are computed by solving the IVP \eqref{eq:conslaw} as described in Definition \ref{def:activefluxtraditional} over half the time step and over the full time step (from the same initial data $q_\text{recon}$), and then these values are used to compute fluxes in the update of the average:
\begin{align}
 \bar q_i^{n+1} &= \bar q_i^n - \Delta t \sum_{\ell = 0}^2 \omega_\ell \frac{f(q^{n + \frac{\ell}{2}}_{i+\frac12}) - f(q^{n + \frac{\ell}{2}}_{i-\frac12})}{\Delta x}  \label{eq:updateaverageleapfrog}
\end{align}
Here, $\omega_\ell = \Big(\frac{1}{6},\frac{2}{3}, \frac{1}{6}\Big)$ are the quadrature weights of Simpson's rule, which is the adequate quadrature for a third order method.

The original publication \cite{vanleer77} considered Active Flux only for linear advection $\del_t q + c \del_x q = 0$, in which case one has for positive $c$
\begin{align}
 q_{i+\frac12}^{n+\frac{\ell}{2}} &= q^n_{\text{recon},i}\left(x_{i+\frac12} - c \frac{\ell}{2} \Delta t\right) \qquad \ell = 1,2 \label{eq:updatelinearadvectionvanleer}\\
 \bar q_i^{n+1} &= \bar q_i^n - \Delta t \sum_{\ell = 0}^2 \omega_\ell c\frac{q^{n + \frac{\ell}{2}}_{i+\frac12} - q^{n + \frac{\ell}{2}}_{i-\frac12}}{\Delta x}   \label{eq:updateaverageleapfrogvanleer}
\end{align}
Given the unique parabolic reconstruction 
\begin{align}
q^n_{\text{recon},i}(x) = \frac{6 \bar q_i - q_{i-\frac12} - q_{i+\frac12}}{4} + \frac{q_{i+\frac12} - q_{i-\frac12}}{\Delta x}x + 3 \frac{q_{i-\frac12} + q_{i+\frac12} - 2 \bar q_i}{\Delta x^2} x^2 \label{eq:parabolicrecon}
\end{align}
which fulfills
\begin{align}
 q^n_{\text{recon},i}\left( \pm \frac{\Delta x}{2}  \right ) &= q^n_{i\pm\frac12} & 
 \frac{1}{\Delta x} \int_{-\frac{\Delta x}{2}}^{\frac{\Delta x}{2}} \dd x \, q^n_{\text{recon},i}(x)  &= \bar q_i^n
\end{align}
Equations \eqref{eq:updatelinearadvectionvanleer}--\eqref{eq:updateaverageleapfrogvanleer} read
\begin{align}
 q^{n+\frac{\ell}{2}}_{i+\frac12} &= q_{i+\frac12}^n \left( 1-2\nu \ell + \frac34 \nu^2 \ell^2  \right )+q^n_{i-\frac12} \left( -\nu \ell + \frac34 \nu^2 \ell^2  \right ) -\frac32 \nu \ell \bar q^n_i (\nu \ell - 2)\\
 \bar q_i^{n+1} &=  \bar q_i^n  (1 + 2 \nu) (1-\nu)^2 + \bar q_{i-1}^n \nu^2 (3-2\nu)  - q_{i-\frac32}^n (1-\nu) \nu^2 +   q_{i-\frac12}^n  \nu (1-\nu)  -  q_{i+\frac12}^n (1-\nu)^2 \nu
 %&= \nu \ell (2-\nu\ell) q^n_i + \frac12 \nu \ell (\nu \ell-1) q^n_{i-\frac12} + \frac12 (\nu\ell-2)(\nu \ell-1)q^n_{i+\frac12}
\end{align}
with $\Delta t = \nu  \Delta x/c$.% and having in the last step introduced the point value at center via $\bar q_i = \frac{q_{i+\frac12} + 4q_i + q_{i-\frac12}}{6}$. Thus,
%\begin{align}
% \frac{q^{n+\frac{\ell}{2}}_{i+\frac12} - q^n_{i+\frac12}}{\Delta t} &= c \frac{\frac12\ell (2-\nu\ell) q^n_i + \frac14 \ell (\nu \ell-1) q^n_{i-\frac12} + \frac14 \ell (\nu\ell-3)q^n_{i+\frac12}}{\Delta x/2} \\
% &= -c \ell \frac{  \frac14  q^n_{i-\frac12}- q^n_i + \frac34  q^n_{i+\frac12}}{\Delta x/2} + c \nu\ell^2 \frac{\frac14  q^n_{i-\frac12} -\frac12 q^n_i  + \frac14 q^n_{i+\frac12}}{\Delta x/2}\\
% &\simeq - c \frac{\ell}{2} \del_x q + \frac12 c \nu \Delta x \frac{\ell^2}4 \del_x^2 q + \mathcal O(\Delta x^2)
%\end{align}

This method is stable up to $\nu = 1$. The main ingredient is an (approximate) evolution operator, which depends on the PDE at hand. Such operators have been found for nonlinear scalar conservation laws and for hyperbolic systems of conservation laws in one dimension by means of an estimate of the characteristic directions in \cite{barsukow19activeflux}. The main difficulty is that local linearization is not sufficiently accurate to obtain third-order accuracy. For systems in multiple spatial dimensions, the concept of characteristics (along which some quantity is constant) is generally to be replaced by that of characteristic cones, which do not allow for simple transformations to characteristic variables. For linear acoustics, an exact evolution operator has been obtained in \cite{eymann13,barsukow17}. It has been then applied to Active Flux on Cartesian two-dimensional grids in \cite{barsukow18activeflux}, and structure preservation properties of the resulting method have been demonstrated: The Active Flux method has been found to be vorticity preserving / low Mach compliant. Despite promising efforts towards approximate evolution operators for the multi-dimensional Euler equations in e.g. \cite{fan17}, further research is required for hyperbolic systems in multiple spatial dimensions that is beyond the scope of the present paper.

\subsection{Semi-discrete Active Flux} \label{ssec:semidiscretegeneral}

In \cite{abgrall20}, a closely related definition was given, which neither relies on a reconstruction, nor on an evolution operator:

\begin{definition}[Semidiscrete Active Flux from \cite{abgrall20}] \label{def:activefluxabgrall}
The \emph{semidiscrete Active Flux method from \cite{abgrall20}} is the following semi-discretization
\begin{align}
 \begin{cases} \displaystyle \frac{\dd}{\dd t} \bar q_i(t) = - \frac{f(q_{i+\frac12}(t)) - f(q_{i-\frac12}(t))}{\Delta x} \\ \\
 \displaystyle \frac{\dd}{\dd t} q_{i+\frac12}(t) = -R\Big(q_{i+\frac12-k}(t), \bar q_{i-k+1}(t), q_{i+\frac12-k+1}(t), \ldots, \bar q_{i + m}(t), q_{i+\frac12+m}(t)\Big)  \qquad i \in \mathbb Z, k\geq 0, m \geq 0 \end{cases} \label{eq:activefluxdefabgrall}
\end{align}
of \eqref{eq:conslaw} with the interpretations \eqref{eq:traditionaldofinterpretations} and $R$ a consistent approximation of $\del_x f(q)$ at $x_{i+\frac12}$.
\end{definition}

Its order of accuracy depends on the approximation order of $R$ (and eventually on the order of accuracy of the time integration). In \cite{abgrall20}, $R$ is found as follows. First, a point value $q_{i,\text{midpoint}}(t)$ at the center of the cell is introduced and the average expressed as $\bar q_i = \frac{1}{6} (q_{i-\frac12} + 4 q_{i,\text{midpoint}} + q_{i+\frac12})$. This allows to choose a standard finite difference formula from \cite{iserles82} on a grid of half the spacing:
\begin{align}
D_3(q_{i-\frac12}, q_{i,\text{midpoint}}, q_{i+\frac12},q_{i+1,\text{midpoint}}) := \frac{\frac16 q_{i-\frac12} - q_{i,\text{midpoint}} +\frac12 q_{i+\frac12} + \frac13 q_{i+1,\text{midpoint}}}{\Delta x/2} \simeq \del_x q + \mathcal O(\Delta x^{3})
\end{align}

For linear advection $\del_t q + c \del_x q = 0$, the equation for the point value at cell interface then reads
\begin{align}
 \frac{\dd}{\dd t} q_{i+\frac12}(t) &= - c \frac{\frac16 q_{i-\frac12}(t) - q_{i,\text{midpoint}}(t) +\frac12 q_{i+\frac12}(t) + \frac13 q_{i+1,\text{midpoint}(t)} }{\Delta x / 2} \label{eq:updateabgrall}
\end{align}
To the same order of accuracy, it can be expressed in the form \eqref{eq:activefluxdefabgrall} with $k=1$, $m = 1$ as
\begin{align}
 R(q_{i+\frac12-k}{t}, \bar q_{i-k+1}, q_{i+\frac12-k+1}(t), \ldots, \bar q_{i + m}, q_{i+\frac12+m}) = c\frac{\frac{5}{6} q_{i-\frac12} -3 \bar q_i + \frac43 q_{i+\frac12}  + \bar q_{i+1} - \frac{1}{6}q_{i+\frac32}  }{\Delta x } \label{eq:iserles3viaavg}
\end{align}

In \cite{abgrall20}, system \eqref{eq:activefluxdefabgrall} is integrated with a Runge-Kutta method. The stability bound for the above method is $\nu_\text{max} =\frac{c \Delta t_\text{max}}{\Delta x} = 0.77$ (see Section \ref{app:stability}). This method has been successfully applied to the Euler equations in \cite{abgrall20} by replacing the prefactor $c$ in \eqref{eq:updateabgrall} by the Jacobian $J = f'$ and using two finite-difference formulae (biased in different directions) in order to include upwinding.
%\begin{align}
%\frac{\dd}{\dd t} q_{i+\frac12}(t) &= - J(q_{i+\frac12}) \frac{\frac16 q_{i-\frac12} - q_i +\frac12 q_{i+\frac12} + \frac13 q_{i+1} }{\Delta x / 2} \\
%\end{align}
%where $J = f'$ is the Jacobian.

Later, (Sections \ref{sec:points} and \ref{sec:moments}) we discuss extensions of Active Flux to higher order that include other types of degrees of freedom. In these cases, the definitions are modified accordingly in order to include the respective updates of new degrees of freedom. This is discussed in the corresponding Sections.

\subsection{Limiting} \label{ssec:limiting}

Limiting is necessary for Active Flux because it is a high order method, and because it is linear when applied to linear problems. The usual philosophy of limiting in the context of finite volume methods is to modify the reconstruction such that its value at the cell interface satisfies a maximum principle. Examples of such approaches are slope limiters, or the approach in \cite{zhang10} where the reconstruction polynomial is scaled around its average. For Active Flux, it is not possible to combine such ways of limiting with continuity of the reconstruction because the values at cell interfaces are prescribed. Limiting strategies for Active Flux that give up continuity can be found in e.g. \cite{kerkmann18,chudzik21}. 

Here, we prefer to maintain continuity. The traditional Active Flux method uses a reconstruction to evolve the point values, and it is a natural approach to limiting to replace an oscillative reconstruction by a monotone, or less oscillative function. Such limiting strategies based on piecewise defined functions (e.g. a constant and a parabola) have been proposed in \cite{roe15,barsukow20swaf}. A simple limiting strategy has been introduced in \cite{barsukow19activeflux} where the parabolic reconstruction is replaced by a power law. This shall be taken as inspiration here. 

Recall the following result from \cite{barsukow19activeflux}:

\begin{theorem}
Consider monotone data
\begin{align}
q_{i-\frac12} \leq \bar q_i \leq q_{i+\frac12} \quad \text{or} \quad q_{i-\frac12} \geq \bar q_i \geq q_{i+\frac12} \label{eq:monotonedata}
\end{align}
Their unique parabolic interpolant is monotone, iff
\begin{align}
  r := \frac{q_{i+\frac12} - \bar q_i }{ \bar q_i - q_{i-\frac12} } \in \left[\frac12, 2\right] \label{eq:monotnicityconditionparabola}
\end{align}
\end{theorem}

\begin{figure}
 \centering
 \includegraphics[width=0.19\textwidth]{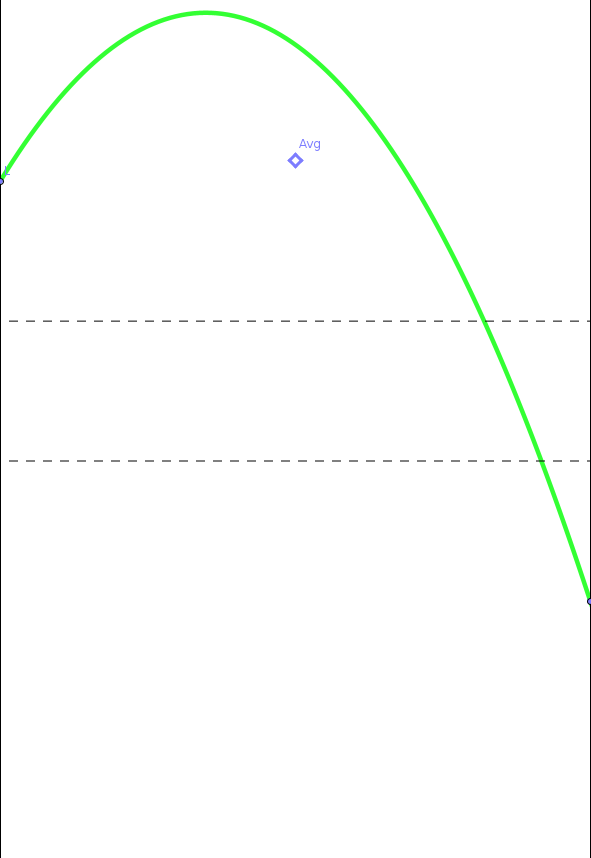}\hfill\includegraphics[width=0.19\textwidth]{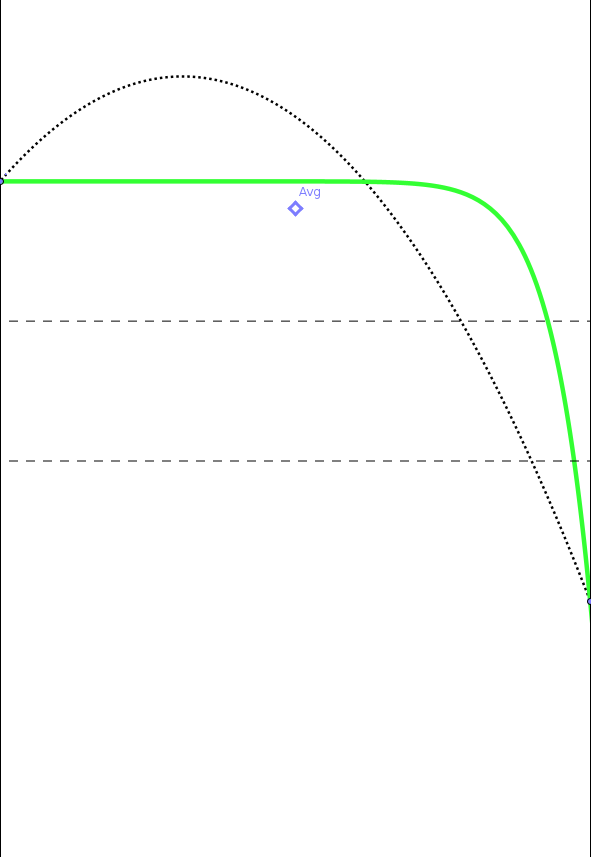}\hfill\includegraphics[width=0.19\textwidth]{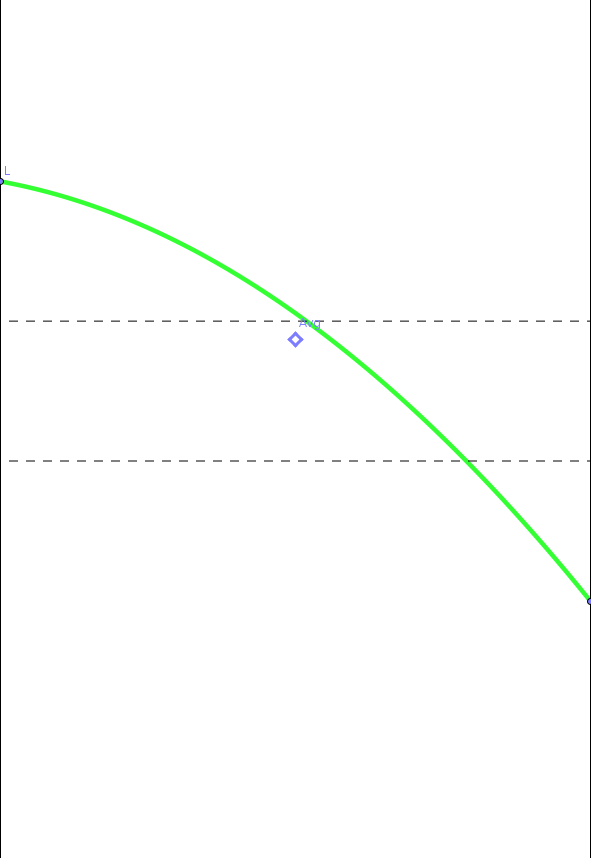}\hfill\includegraphics[width=0.19\textwidth]{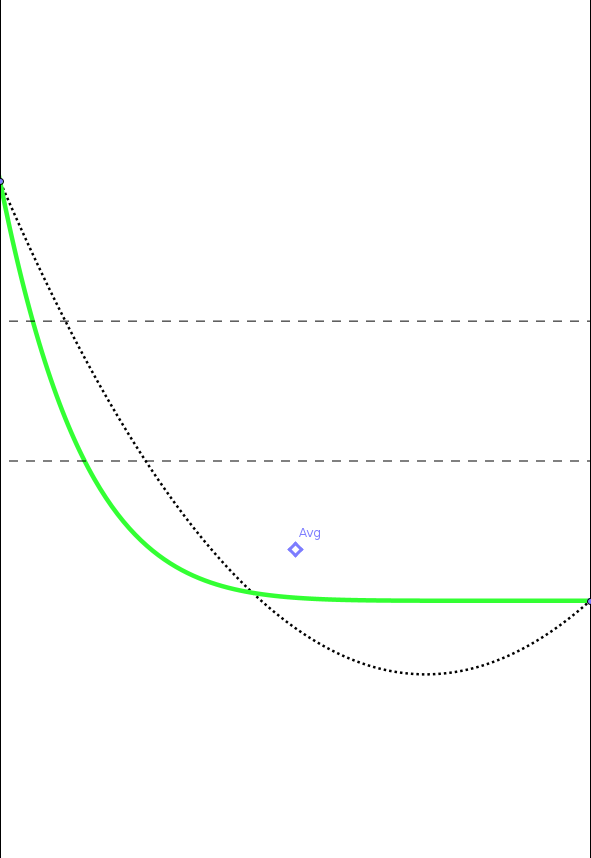}\hfill\includegraphics[width=0.19\textwidth]{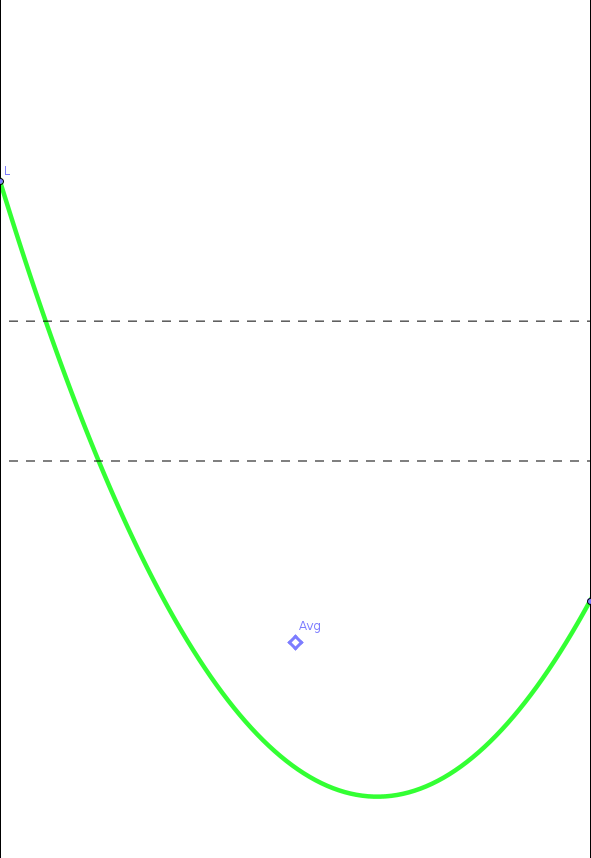}
 \caption{Different cases of limiting. The two point values remain the same, the diamond in the cell center indicates the value of the cell average. The green solid line shows the choice of reconstruction, the dotted line shows the parabolic reconstruction whenever it is not the one chosen. The dashed horizontal lines indicate the region of the cell average inside which the parabolic reconstruction is monotone. \emph{From left to right}: The cell average is larger than the two point values, there exists is no continuous monotone reconstruction, and the parabola is accepted. The cell average is between the two point values, but the parabola is not monotone, we choose one of the power laws. The average is inside the region given by the dashed lines, we choose the (monotone) parabola. The average is below the region, we choose the other power law. The average is below both point values, we choose the parabola again.}
 \label{fig:limitingpowerlaw}
\end{figure}

Thus, parabolic interpolants can violate the maximum principle if the average $\bar q_i$ is close to one of the point values $q_{i\pm\frac12}$. The power law reconstructions, proposed in \cite{barsukow19activeflux}

\begin{align}
 q_{\text{recon},i,\text{power-law},1}(x) &= q_{i-\frac12} + (q_{i+\frac12} - q_{i-\frac12}) \left (   \frac{x + \frac{\Delta x}{2}}{\Delta x} \right )^\frac{q_{i+\frac12} - \bar q_i}{\bar q_i - q_{i-\frac12}} \label{eq:limitingpowerlaw1}\\
 q_{\text{recon},i,\text{power-law},2}(x) &= q_{i+\frac12} - (q_{i+\frac12} - q_{i-\frac12}) \left (   \frac{\frac{\Delta x}{2} - x}{\Delta x} \right )^\frac{\bar q_i - q_{i-\frac12}}{q_{i+\frac12} - \bar q_i}   \label{eq:limitingpowerlaw2}\\
 \nonumber x &\in \left[ - \frac{\Delta x}{2}, \frac{\Delta x}{2}  \right]
 \end{align}
are interpolating the same data while being monotone whenever the data are, and are thus a good replacement for the parabolic reconstruction whenever the latter fails to be monotone. Observe that the expression appearing in \eqref{eq:monotnicityconditionparabola} is just the power in \eqref{eq:limitingpowerlaw1}. Thus, when  
\begin{align}
 r = \frac{q_{i+\frac12} - \bar q_i }{ \bar q_i - q_{i-\frac12} } = 2 \qquad \Rightarrow \qquad \bar q_i - q_{i-\frac12} = \frac{q_{i+\frac12} - q_{i-\frac12}}{3}
\end{align}
the power law \eqref{eq:limitingpowerlaw1} reduces to the parabolic interpolant. The same is true for \eqref{eq:limitingpowerlaw2} for
\begin{align}
 r = \frac{q_{i+\frac12} - \bar q_i }{ \bar q_i - q_{i-\frac12} } = \frac12 \qquad \Rightarrow \qquad \bar q_i - q_{i+\frac12} = \frac{q_{i-\frac12} - q_{i+\frac12}}{3}
\end{align}
If the data are monotone, then $r \geq 0$. For $r \in [0, \frac12)$ one therefore should choose the power-law \eqref{eq:limitingpowerlaw2}, which at $r=\frac12$ becomes the parabola. Then, for $r \in [\frac12, 2]$, the parabola is monotone, and at $r = 2$ it is now the power-law \eqref{eq:limitingpowerlaw1} which is equal to the parabola. For $r \in (2, \infty)$ one thus chooses the power law \eqref{eq:limitingpowerlaw1}:

\begin{theorem}
 Consider monotone data fulfilling \eqref{eq:monotonedata}. The reconstruction
 \begin{align}
  q_{\text{recon},i}(x) &= \begin{cases} q^n_{\text{recon},i,\text{power-law},1}(x) & |\bar q_i - q_{i-\frac12}| < \frac{|q_{i+\frac12} - q_{i-\frac12}|}{3} \\
  q^n_{\text{recon},i,\text{power-law},2}(x) & |\bar q_i - q_{i+\frac12}| < \frac{|q_{i+\frac12} - q_{i-\frac12}|}{3} \\
  \text{parabolic} & \text{else}
  \end{cases} \qquad \nonumber x \in \left[ - \frac{\Delta x}{2}, \frac{\Delta x}{2}  \right]
 \end{align}
 is continuous as a function of $\bar q_i$.
\end{theorem}

For data that are not monotone (i.e. violating \eqref{eq:monotonedata}), a parabolic reconstruction is used. See Figure \ref{fig:limitingpowerlaw} for an overview.

A rational interpolation (suggested in \cite{kerkmann18}) is similar in spirit to the power law reconstruction, but its coefficients cannot be computed analytically because of the condition on its average, which requires integration of a rational function.

Note that in all suggestions of limiting available in the literature, only the update of the point values is limited. This provides satisfactory results, but does not exclude the appearance of oscillations, as the finite volume step 
\begin{align}
 \frac{\dd}{\dd t} \bar q_i(t) = - \frac{f(q_{i+\frac12}(t)) - f(q_{i-\frac12}(t))}{\Delta x}
\end{align}
remains unlimited. A limiting of the finite volume step is subject of future work.

\subsection{Transonic upwinding} \label{ssec:transonicupwindinggeneral}

\begin{figure}
 \centering
 \includegraphics[width=0.3\textwidth]{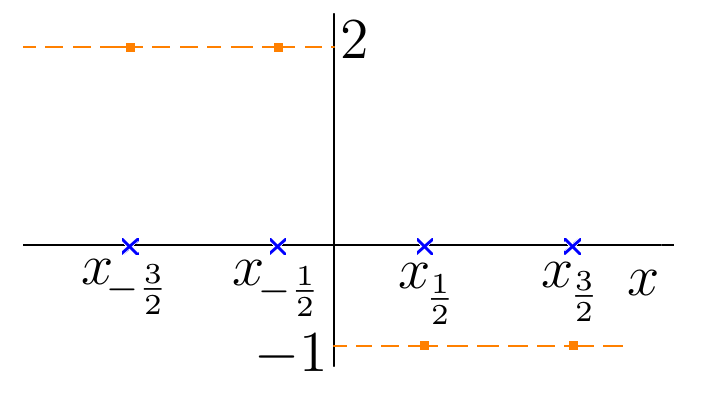}
 \caption{Setup for the illustration of problems if transonic upwinding is not performed carefully.}
 \label{fig:setuptransonicupwinding}
\end{figure}

In \cite{kerkmann18,barsukow19activeflux}, the following problematic situation has been observed with early versions of Active Flux. Consider initial data 
\begin{align}
q_0(x) = \begin{cases}
            2 & x \leq 0  \\ -1 & x > 0
          \end{cases} \label{eq:kerkmanntest}
\end{align}
for Burgers' equation, such that the discrete data initially are (see also Figure \ref{fig:setuptransonicupwinding})
\begin{align}
\ldots = q_{-\frac32} = q_{-\frac12} &= 2 & q_{\frac12} = q_\frac32 = \ldots &= -1\\
\ldots = \bar q_{-2} = \bar q_{-1} &= 2 & \bar q_{1} = \bar q_{2} = \ldots &= -1\\
q_{\text{recon},-1} &\equiv 2 & q_{\text{recon},1} &\equiv -1 
\end{align}
The value of $\bar q_0$ is irrelevant for the discussion and any can be chosen. This value depends on where exactly inside the cell the location $x=0$ of the shock will be. In other words, the problem discussed below is generic and will not disappear if the location of the shock is changed slightly.

The exact solution of \eqref{eq:kerkmanntest} would be a shock moving at speed $\frac12$. In the update of $q_{-\frac12}$, if its value 2 is used to estimate the upwind direction, information would be taken from the left, where the reconstruction is identically equal to 2. $q_{-\frac12}$ would thus remain stationary. The same happens to $q_{\frac12}$, whose initial value is negative, and thus information is taken from the right. In fact, the reconstruction $q_{\text{recon},0}$ thus never contributes. The numerical method keeps all point values stationary, while the average $\bar q_0$ in cell 0 will grow, because the flux difference is $(q_{\frac12})^2/2 - (q_{-\frac12})^2/2 \neq 0$.

\begin{figure}
 \centering
 \includegraphics[width=0.7\textwidth]{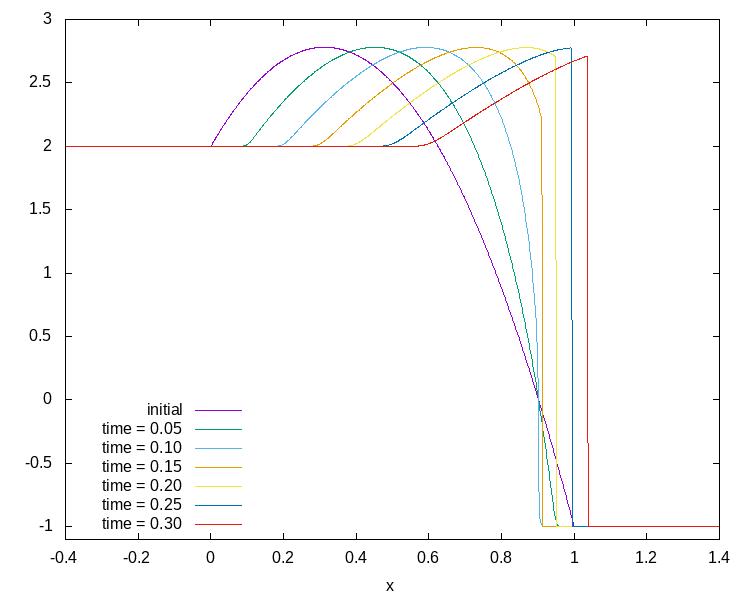}
 \caption{The evolution of continuous initial data corresponding to a parabolic reconstruction of \eqref{eq:kerkmanntest} with $\Delta x = 1$ for Burgers' equation. The time step based on the maximum of the reconstruction is $\Delta t \simeq 0.357$, and even larger if based on the point values and averages alone. At $x = 1$ one observes initially an information flow from the right, but as the shock forms, it stalls and starts moving to the right. After $t \simeq 0.25$, the information is flowing from the left at $x = 1$.}
 \label{fig:burgers-kerkmann-shock}
\end{figure}

The error is the wrong estimate of the upwind direction at $x_{\frac12}$. An exact time evolution of the continuous reconstruction $q_{\text{recon}}$ would see it form a shock, and this shock would first move to the left, and then to the right, reaching $x_{\frac12}$ before the end of the time step. At this moment, the direction of information flow at $x_{\frac12}$ would be reversed, and information would no longer be flowing from the right (see Figure \ref{fig:burgers-kerkmann-shock}). Thus, the deficiency of the early versions of Active Flux is due to an erroneous estimate of the upwind direction in the point value update. This estimate needs to take into account the fact that a continuous reconstruction might self-steepen before the end of the time step.

This error is most visible in the transonic case, because the continuous reconstruction has time to self-steepen during the time step. Depending on the CFL condition, setups close to transonic might be affected as well. In the worst case, the consequence is a shock that is not moving (and thus violating the Rankine-Hugoniot conditions), with the cell average in its vicinity growing without bound. In setups that are not strictly transonic, an artificially large average at the location of the shock is visible in the form of a spike (see examples in \cite{barsukow19activeflux}).

In \cite{barsukow19activeflux} it has been suggested to estimate the upwind direction more carefully. In particular, when approximating iteratively the foot point of the characteristic, it has been suggested to initialize the iteration with not only the value at the respective cell interface, but also with values from its neighbours. If a shock forms over the time step, one would thus obtain several candidate characteristics, while in a smooth setup the iteration would converge to the unique characteristic. A simple strategy that works well in practice is then to select from the candidate characteristics the one whose speed is largest in modulus. In \cite{barsukow19activeflux,barsukow20swaf} it has been demonstrated that this strategy works well even for systems. Some authors (e.g. \cite{chudzik21}) suggest to average to values obtained from the candidate characteristics.

\section{High order via larger stencils} \label{sec:findiff}

\subsection{Discretization in space}

The first strategy for high order accuracy that we propose involves the traditional degrees of freedom (cell averages and shared point values at cell interfaces), and uses larger stencils. The starting point therefore is the semidiscrete Active Flux of Definition \ref{def:activefluxabgrall}, with a Runge-Kutta time integration. Observe again that the update of the averages 
\begin{align}
 \frac{\dd}{\dd t} \bar q_i(t) = - \frac{f(q_{i+\frac12}(t)) - f(q_{i-\frac12}(t))}{\Delta x} \label{eq:averageupdatefindiff}
\end{align}
is exact, and that the order of accuracy of the method is given by the order of accuracy of the right-hand side of the point value update in \eqref{eq:activefluxdefabgrall} and of the order of accuracy of the time integration.

Consider linear finite difference approximations to the derivative of the form
% &\multicolumn{2}{c}{a_\frac32}
\begin{align}
 \text{FD}p = \frac{1}{\Delta x} \sum_{j = - \ell}^\ell (b_j^{(0)}  \bar q_{i+j} + b_{j+\frac12} q_{i + j+\frac12}) = q'( x_{i+\frac12}) + \mathcal O(\Delta x^{p-1}) \label{eq:finitedifferencegeneral}
\end{align}
The definition of $p$ is such that a finite difference formula FD$p$ will give rise to a $p$-th order Active Flux method, because it is sufficient to use a derivative approximation that is exact on polynomials of degree $p -1$.

We shall adopt the following tableau notation for such a finite difference formula that involves point values at cell interfaces and cell averages

\begin{notation}

\begin{align}
 \frac{1}{\Delta x} \sum_{j = - \ell}^\ell (b_j^{(0)} \bar q_{i+j} + b_{j+\frac12} q_{i + j+\frac12}) \equiv \phantom{mmmmmmmmmmmmmmmmmmmmmmmmm}\label{eq:definitiontableau} \\\nonumber
 \begin{array}{|cc|ccc|ccc|ccc||ccc|ccc|ccc|c|}
    \hline
    b_{-\ell}^{(0)} &\phantom{s}&&                     \cdots &\phantom{i}&&                        b_{-1}^{(0)} &\phantom{i}&&                        b_0^{(0)} &&&                      b_1^{(0)} &&&   \cdots&&&   b_\ell^{(0)}&\phantom{i}&\\
    &\multicolumn{2}{c}{b_{-\ell+\frac12}}&  &\multicolumn{2}{c}{b_{-\frac32}}&      &\multicolumn{2}{c}{b_{-\frac12}}&      &\multicolumn{2}{c}{b_\frac12}&    &\multicolumn{2}{c}{b_\frac32}&    &\multicolumn{2}{c}{b_{\ell-\frac12}}&  &\multicolumn{2}{c|}{b_{\ell + \frac12}} \\[4pt]\hline
 \end{array} 
\end{align}
The vertical lines indicate cell interfaces (such that coefficients of cell averages are written ``inside'' the cell and coefficients of point values are associated to cell interfaces). The double vertical line indicates the cell interface at which the finite difference is supposed to provide an approximation to the derivative (of the corresponding order of accuracy). Coefficients not marked explicitly in the tableau are assumed zero. The notation will be slightly expanded in Section \ref{sec:moments}, where the superscript $^{(0)}$ on the coefficients of the cell averages will become clear.

Here are some examples ($a_2, a_4, a_6$ are free parameters):

\begin{align}
 \text{FD}2 &= \begin{array}{|c|ccc||c|}
    \hline &&&&\\[-12pt]
    \phantom{mn}&& 2-2a_2 && \\
    \multicolumn{2}{|c}{-2 + a_2} && \multicolumn{2}{c|}{a_2}\\[4pt]\hline
 \end{array}\label{eq:stencil1st} \\
 \text{FD}4a &=\begin{array}{|c|ccc||ccc|c|}
  \hline &&&& &&& \\[-12pt] 
  \phantom{m}&& - 2 - \frac{3a_4}{4}  &&& 2-\frac{3a_4}{4}  &\phantom{n}& \\
  \multicolumn{2}{|c}{\frac{2+a_4}{4}} && \multicolumn{2}{c}{a_4} && \multicolumn{2}{c|}{\frac{-2 + a_4}{4}} \\[4pt]\hline
 \end{array}\label{eq:stencil3rd}\\
 \text{FD}6b &=\begin{array}{|c|ccc|ccc||ccc|c|} 
               \hline &&&& &&& &&&\\[-12pt]
               \phantom{m}&&   \frac{19-22a_6}{54}  &&&   -\frac{89 + 76a_6}{54}  &&&  \frac{50 -11a_6}{27}  &\phantom{m}& \\
               \multicolumn{2}{|c}{\frac{a_6-1}{9}} && \multicolumn{2}{c}{a_6 } && \multicolumn{2}{c}{a_6 } && \multicolumn{2}{c|}{\frac{a_6-4}{9} } \\[4pt] \hline
              \end{array}\label{eq:stencil5th}
              %
%  \text{FD}8a &= \begin{array}{|c|ccc|ccc|ccc||ccc|c|} 
%               \hline &&&& &&& &&& &&&\\[-12pt]
%\phantom{n}&&  \frac{49}{72} - \frac{25 a_7}{96}  &\phantom{mi}&& \frac{293}{72} - \frac{185 a_7}{96} &\phantom{m}&& - \frac{31}{72} - \frac{185 a_7}{96} &&& \frac{436-75a_7}{288} &\phantom{m}& \\\multicolumn{2}{|c}{\frac{-8+3 a_7}{48}} && \multicolumn{2}{c}{-\frac{7}{3} + a_7} && \multicolumn{2}{c}{-3 +\frac{9a_7}{4}} && \multicolumn{2}{c}{a_7 } && \multicolumn{2}{c|}{-\frac13 + \frac{a_7}{16} }  \\[4pt] \hline
%  \end{array}\label{eq:stencil7th}
 \end{align}
 
Note once more that the double vertical line indicates the cell interface at which the finite difference is providing the desired approximation of the derivative. E.g. for FD2,
\begin{align}
 q'(x_{i+\frac12}) &\simeq (-2 + a_2) q_{i-\frac12}  + (2-2a_2) \bar q_i + a_2 q_{i+\frac12}\\
\intertext{while for FD4a,}
 q'(x_{i+\frac12}) &\simeq \frac{2+a_4}{4} q_{i-\frac12}  + \left(-2-\frac{3a_4}{4}\right) \bar q_i + a_4 q_{i+\frac12} + \left(2-\frac{3a_4}{4}\right) \bar q_{i+1} + \frac{-2+a_4}{4} q_{i+\frac32}
\end{align}

\end{notation}

Most finite difference formulas that are considered here involve one degree of freedom more than it would be necessary to obtain the desired order of accuracy; the remaining free parameter is used later to optimize stability. These finite difference formulas are different from standard finite differences because they involve point values and averages. Table \ref{tab:finitedifferences} summarizes some of the possible choices. Given a finite difference formula

\newcommand{\findiffextraspace}{\phantom{m}}
\begin{align}
 D = \begin{array}{|cc|ccc|ccc|ccc||ccc|ccc|ccc|c|}
    \hline
    \findiffextraspace\findiffextraspace b_{-\ell}^{(0)} \findiffextraspace &\phantom{s}&&                     \cdots &\phantom{i}&&                        \findiffextraspace b_{-1}^{(0)} \findiffextraspace &\phantom{i}&&                       \findiffextraspace b_0^{(0)}\findiffextraspace &&&                      b_1^{(0)} &&&   \cdots&&&   b_\ell^{(0)}&\phantom{i}&\\
    &\multicolumn{2}{c}{b_{-\ell+\frac12}}&  &\multicolumn{2}{c}{b_{-\frac32}}&      &\multicolumn{2}{c}{b_{-\frac12}}&      &\multicolumn{2}{c}{b_\frac12}&    &\multicolumn{2}{c}{b_\frac32}&    &\multicolumn{2}{c}{b_{\ell-\frac12}}&  &\multicolumn{2}{c|}{b_{\ell + \frac12} \findiffextraspace}  \\[4pt]\hline
 \end{array}
\end{align}
define the flipped formula
\begin{align}
 D^* = \begin{array}{|c|ccc|ccc|ccc||ccc|ccc|ccc|c|}
    \hline
    &&        -b_\ell^{(0)} &\phantom{i}&&                      \cdots &\phantom{i}&&                        -b_1^{(0)} &&&                      -b_0^{(0)} &&&   \cdots&&&   -b_{-\ell}^{(0)}&\phantom{i}&\\
    \multicolumn{2}{|c}{-b_{\ell+\frac12}}&  &\multicolumn{2}{c}{-b_{\ell - \frac12}}&      &\multicolumn{2}{c}{-b_{\frac32}}&      &\multicolumn{2}{c}{-b_\frac12}&    &\multicolumn{2}{c}{-b_{-\frac12}}&    &\multicolumn{2}{c}{-b_{-\ell+\frac32}}&  &\multicolumn{2}{c|}{-b_{-\ell + \frac12}} \\[4pt]\hline
 \end{array}
\end{align}
It is an approximation of the derivative at the same location and of the same order of accuracy as $D$, and arises because the spatial derivative changes sign upon reflection $x \mapsto -x$.

We propose to use these finite differences in the update of the point values by writing
\begin{align}
 \frac{\dd}{\dd t} q_{i+\frac12}(t) = - \Big ( f'(\tilde q_{i+\frac12})^+ D  + f'(\tilde q_{i+\frac12})^- D^* \Big ) \label{eq:pointvalueupdatefindiff}
\end{align}
where $D$ is any finite difference formula such as \eqref{eq:stencil1st}--\eqref{eq:stencil5th}. In the scalar case the positive and negative parts are given by
\begin{align}
 f'(\tilde q_{i+\frac12})^+ &= \max(0, f'(\tilde q_{i+\frac12})) &f'(\tilde q_{i+\frac12})^- &= \min(0, f'(\tilde q_{i+\frac12}))
\end{align}

In the case of systems, $f'$ is the Jacobian and they are defined via its eigenvalues:
\begin{align}
 f'(\tilde q_{i+\frac12}) &= R \diag(\lambda_1, \ldots, \lambda_m) R^{-1} & f'(\tilde q_{i+\frac12})^\pm &:= R \diag(\lambda_1^\pm, \ldots, \lambda_m^\pm) R^{-1}
\end{align}

It remains to define $\tilde q_{i+\frac12}$. For scalar conservation laws ($m=1$), if $f'$ does not switch sign in the computational domain, then the natural choice $\tilde q_{i+\frac12} =  q^n_{i+\frac12}$ works. Further details are given in the next section.

The use of finite differences in the point value update is inspired by the approach in \cite{abgrall20}, but different in the following crucial aspect. In \cite{abgrall20}, a point value $q_{i,\text{midpoint}}$ at the cell center is first estimated through Simpson's rule:
\begin{align}
 \bar q_i =: \frac{q_{i-\frac12} + 4 q_{i,\text{midpoint}} + q_{i+\frac12}  }{6} \label{eq:tmpsimpsonpointvaueatcellcenter}
\end{align}
One thus obtains a grid of half spacing, on which standard finite differences can be used. (This is described in Section \ref{ssec:semidiscretegeneral}.) While for third order, Simpson's rule is accurate enough, for higher orders the natural expression of the cell average in terms of point values would require solving a linear system, i.e. it would become nonlocal. Indeed, a formula such as
\begin{align}
 %\bar q_i =: \frac{q_{i-\frac32} + 4 q_{i-1,\text{midpoint}} + q_{i-\frac12} + 4 q_{i,\text{midpoint}} + q_{i+\frac12}  }{6}
 \bar q_i =: \frac{- 4 q_{i-1,\text{midpoint}} + 34 q_{i-\frac12} + 144 q_{i,\text{midpoint}} + 34 q_{i+\frac12} - q_{i+1,\text{midpoint}} }{180} = \frac{1}{\Delta x} \int_{-\frac{\Delta x}{2}}^{\frac{\Delta x}{2}} \dd x \, q + \mathcal O(\Delta x^6)
\end{align}
cannot be as easily solved for the midpoint values as \eqref{eq:tmpsimpsonpointvaueatcellcenter}. One could use a quadrature for the average that involves precisely one midpoint value and correspondingly more of the point values at cell interfaces:
\begin{align}
 %\bar q_i =: \frac{q_{i-\frac32} + 4 q_{i-1,\text{midpoint}} + q_{i-\frac12} + 4 q_{i,\text{midpoint}} + q_{i+\frac12}  }{6}
 \bar q_i =: \frac{- 4 q_{i-\frac32} + 189 q_{i-\frac12} + 704 q_{i,\text{midpoint}} + 189 q_{i+\frac12} - q_{i+\frac32} }{1080} = \frac{1}{\Delta x} \int_{-\frac{\Delta x}{2}}^{\frac{\Delta x}{2}} \dd x \, q + \mathcal O(\Delta x^6)
\end{align}
It is clear, however, that the form \eqref{eq:finitedifferencegeneral} includes this case, and is more general.

\subsection{Time integration and stability conditions}

\begin{figure}
 \centering
\begin{minipage}{.5\textwidth}
  \centering
 \includegraphics[width=.8\textwidth]{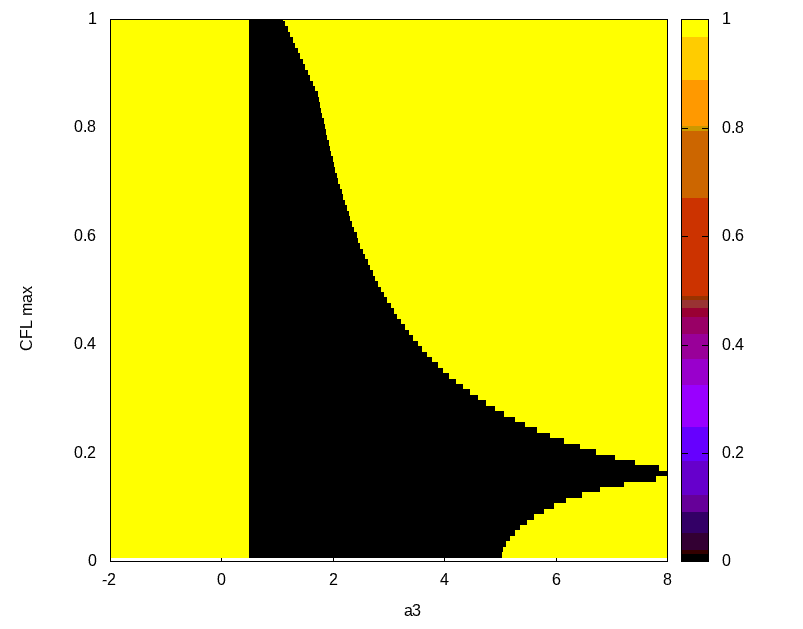} 
 \captionof{figure}{Stability region (RK3) of finite difference formula FD4b (black $=$ stable).}
 \label{fig:rk3stabilityregions3upwind1}
\end{minipage}%
\begin{minipage}{.5\textwidth}
  \centering
 \includegraphics[width=.8\textwidth]{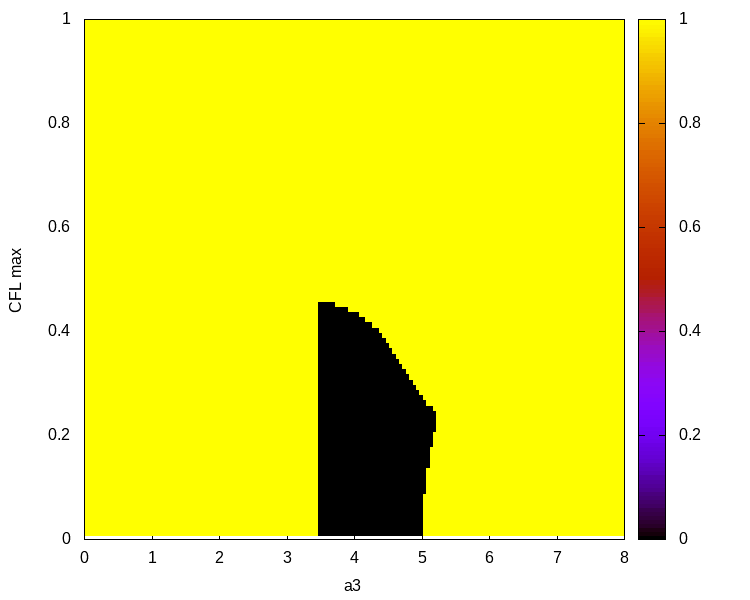} 
 \captionof{figure}{Stability regions (RK3) of finite difference formula FD4c (black $=$ stable).}
 \label{fig:rk3stabilityregions3upwind2}
\end{minipage}
\end{figure}

\begin{figure}  % iserles4-fd.dat / iserles4-fd-upwind1.dat
 \centering
\begin{minipage}{.5\textwidth}
  \centering
 \includegraphics[width=.8\textwidth]{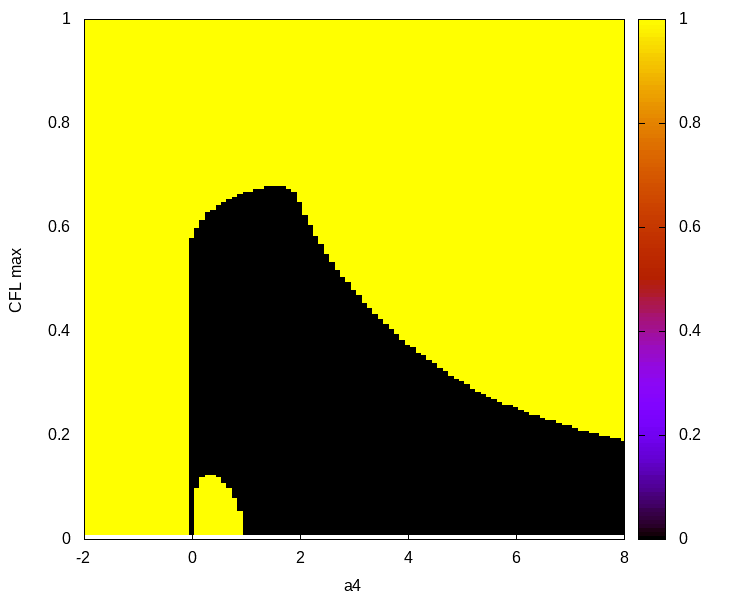} 
 \captionof{figure}{Stability region (RK3) of finite difference formula FD5a (black $=$ stable).}
 \label{fig:rk3stabilityregions4}
\end{minipage}%
\begin{minipage}{.5\textwidth}
  \centering
 \includegraphics[width=.8\textwidth]{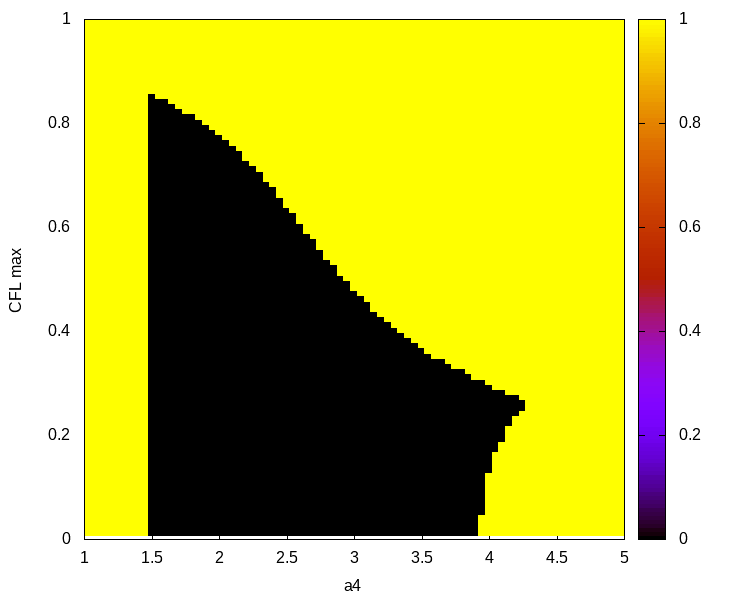} 
 \captionof{figure}{Stability regions (RK3) of finite difference formula FD5b (black $=$ stable).}
 \label{fig:rk3stabilityregions4upwind1}
\end{minipage}
\end{figure}

\begin{figure}
 \centering
\begin{minipage}{.5\textwidth}
  \centering
 \includegraphics[width=0.8\textwidth]{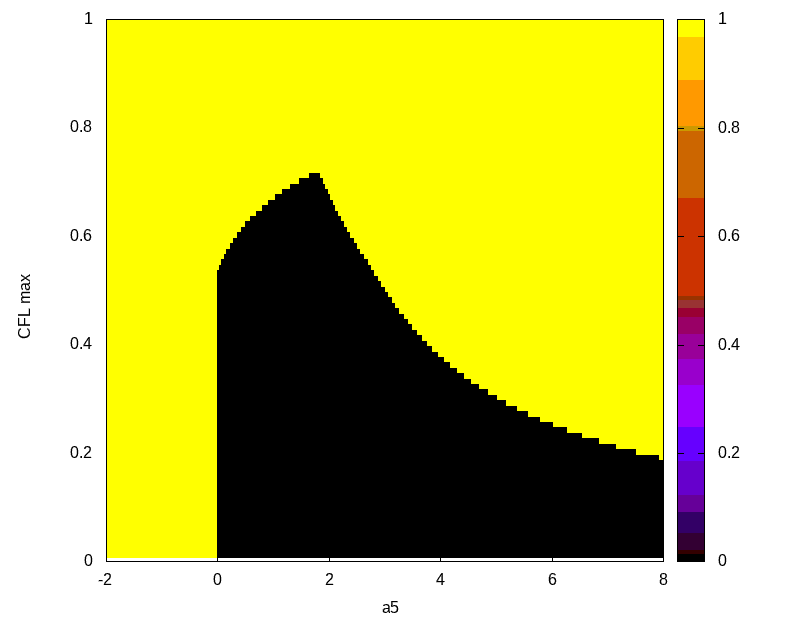} 
 \captionof{figure}{Stability region (RK3) of finite difference formula FD6a (black $=$ stable).}
 \label{fig:rk3stabilityregions5thupwindminus1}
\end{minipage}%
\begin{minipage}{.5\textwidth}
  \centering
 \includegraphics[width=0.8\textwidth]{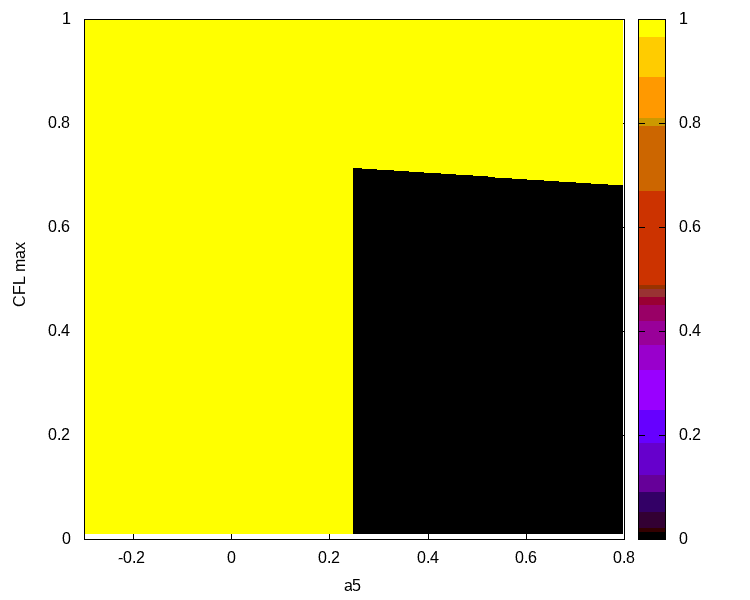}  
 \captionof{figure}{Stability region (RK3) of \eqref{eq:stencil5th} (FD6b) (black $=$ stable).}
 \label{fig:rk3stabilityregions5th}
\end{minipage}
\end{figure}

\begin{figure}
 \centering
 \begin{minipage}{.5\textwidth}
  \centering
 \includegraphics[width=0.8\textwidth]{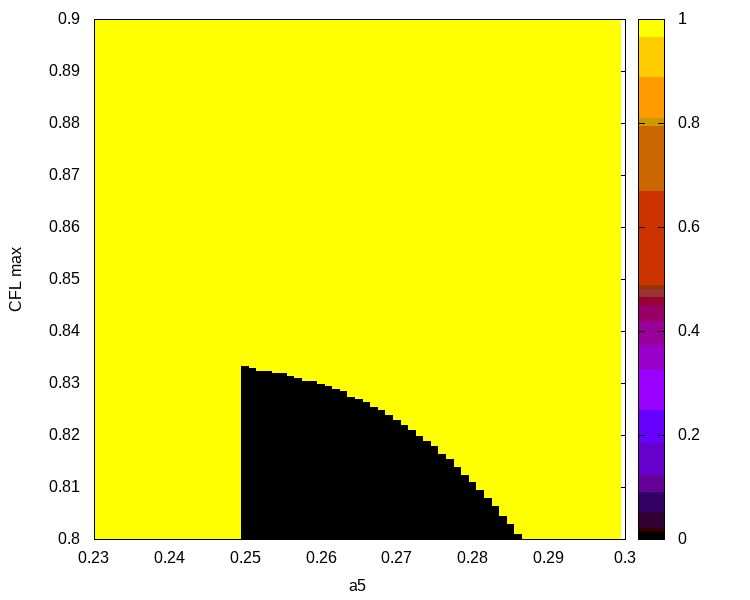} 
 \captionof{figure}{Stability region of \eqref{eq:stencil5th} (FD6b) for RK5 (black $=$ stable).}
 \label{fig:rk5stabilityregions5th}
\end{minipage}%
\begin{minipage}{.5\textwidth}
  \centering
 \includegraphics[width=0.8\textwidth]{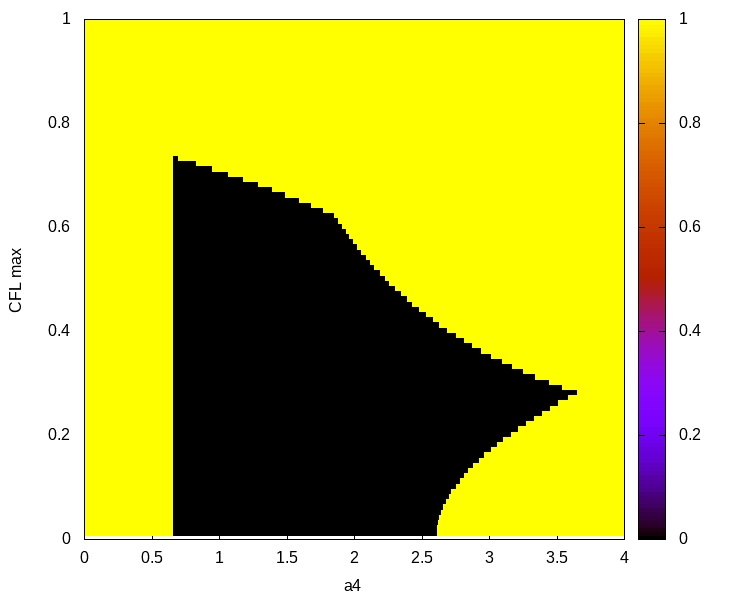}
 \captionof{figure}{Stability region (RK3) of finite difference formula FD7 (black $=$ stable).}
 \label{fig:rk3stabilityregions6th}
\end{minipage}%
\end{figure}

\begin{figure}
 \centering
 \begin{minipage}{.5\textwidth}
  \centering
 \includegraphics[width=0.8\textwidth]{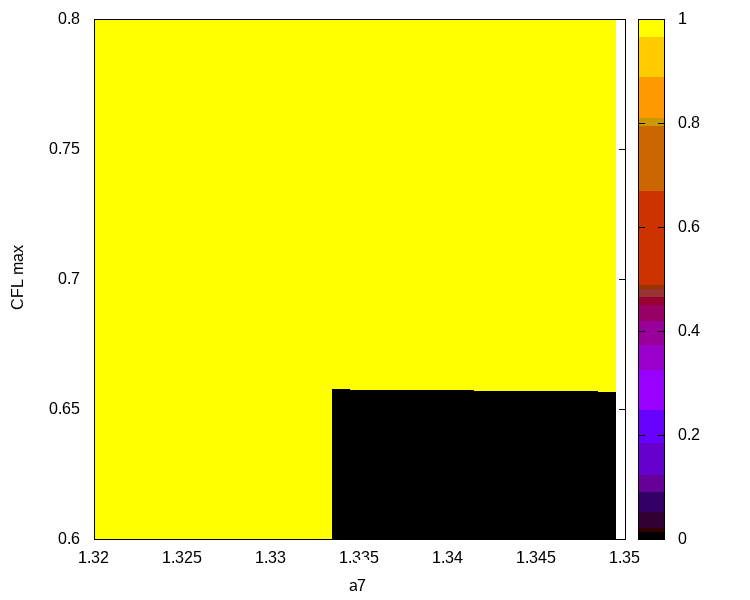}
 \captionof{figure}{Stability region (RK3) of finite difference formula FD8a (black $=$ stable).}
 \label{fig:rk3stabilityregions7th}
\end{minipage}%
\begin{minipage}{.5\textwidth}
  \centering
 \includegraphics[width=0.8\textwidth]{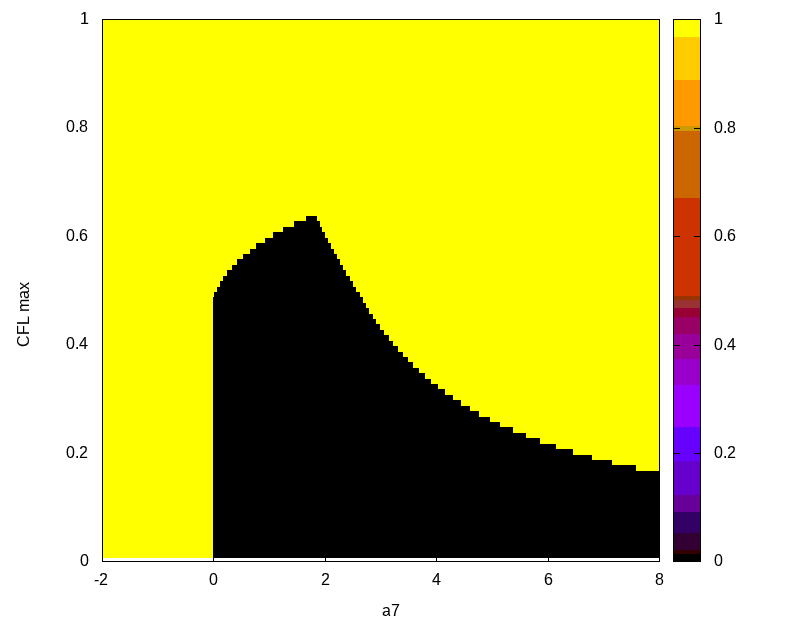} 
 \captionof{figure}{Stability region (RK3) of finite difference formula FD8c (black $=$ stable).}
 \label{fig:rk3stabilityregions7thupwindminus1}
\end{minipage}
\end{figure}

The algorithm consists of the ODE system \eqref{eq:averageupdatefindiff} and \eqref{eq:pointvalueupdatefindiff}. We propose to solve it with a Runge-Kutta time integrator, e.g. a strong stability preserving method of order 3 (SSP-RK3). For linear problems, of course, higher order Runge-Kutta methods are easily available. 
(Compare Figures \ref{fig:rk3stabilityregions5th} and \ref{fig:rk5stabilityregions5th}, where stability of FD6b is shown for RK3 and RK5, respectively, with the maximum CFL number being only slightly higher for RK5.). 
For nonlinear problems, often a lower order time integrator is used with a time step small enough, such that the global error is dominated by the spatial accuracy. We thus content ourselves with results on RK3. The finite difference formulas considered above are typically one cell larger than it would be necessary to obtain the desired order of accuracy, such that the free parameter can be used to optimize stability. The results of the linear stability analysis (see Section \ref{app:stability}) can thus be depicted in a diagram whose axes show the CFL number and the value of the free parameter. Some of these diagrams are shown in Figures \ref{fig:rk3stabilityregions3upwind1}--\ref{fig:rk3stabilityregions7thupwindminus1}. Table \ref{tab:finitedifferences}, shows the maximum CFL numbers for parameter values optimizing stability. For nonlinear problems, the time step is based on the largest characteristic speed $\lambda_\text{max}$ in absolute value, evaluated at the point values. The time step is then given by $\Delta t = \text{CFL} \Delta x / \lambda_\text{max}$.

\begin{figure}
 \centering
\begin{minipage}{.5\textwidth}
  \centering
 \includegraphics[width=0.8\textwidth]{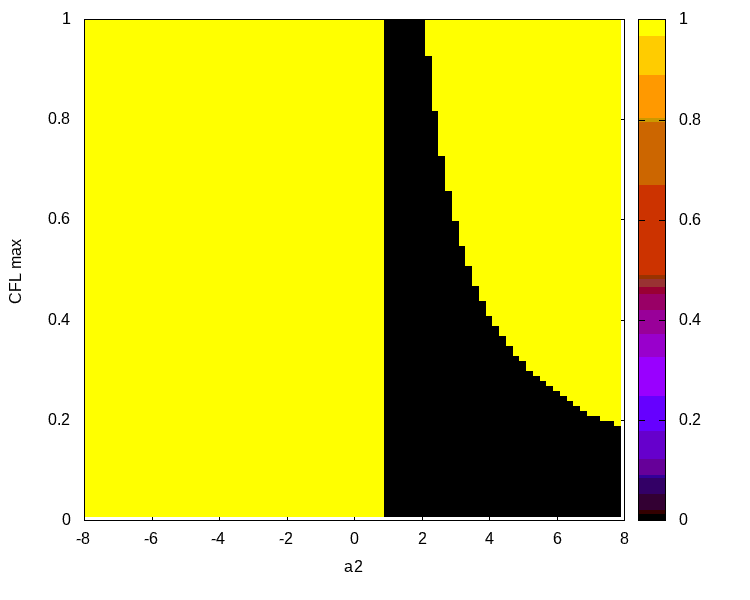} 
 \captionof{figure}{Stability region (RK3) of \eqref{eq:stencil1st} (FD2) (black $=$ stable).}
 \label{fig:rk3stabilityregions1st}
\end{minipage}%
\begin{minipage}{.5\textwidth}
  \centering
 \includegraphics[width=0.8\textwidth]{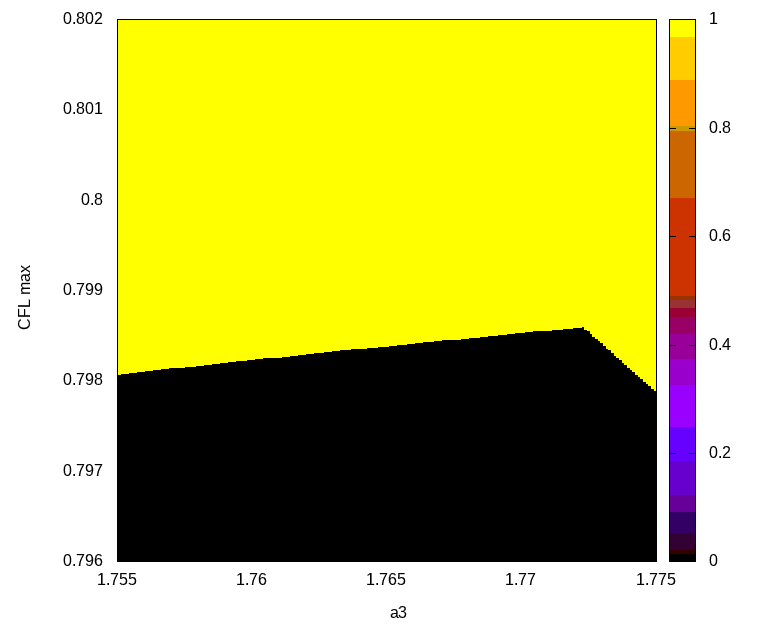}
 \captionof{figure}{Stability region (RK3) of \eqref{eq:stencil3rd} (FD4a) (black $=$ stable).}
 \label{fig:rk3stabilityregions3rd}
\end{minipage}
\end{figure}

\begin{table}
 \centering
\begin{tabular}{c||c||c|c|c||c|c|}
 name           &  FD2/FD3   &    FD4a            & FD4b                     & FD4c            &FD5a       & FD5b\\
 $p$            &  2 (3) &      4             &             4            &      4          & 5         & 5 \\
 \hline\hline
 $b^{(0)}_{-3}$ &        &                    &                          &                 & & \\\hline
 $b_{-\frac52}$ &        &                    &                          &                 & & \\\hline
 $b^{(0)}_{-2}$ &        &                    &                          &                 & & \\\hline
 $b_{-\frac32}$ &        &                    &                          & $a_4-5$         & &$\frac{a_5-3}{3}$ \\\hline
 $b^{(0)}_{-1}$ &        &                    &$\frac16(2-a_4)$          &$\frac{29}2-3a_4$& $-\frac{a_5}{18}$ &$\frac{19}6-\frac{10a_5}9$ \\\hline
 $b_{-\frac12}$ &$-2+a_2$&$\frac{2+a_4}{4}$   & $ a_4-1$                 & $ 4 (a_4 - 4) $ & $\frac{1+a_5}{2}$ &$2 (a_5-2)$      \\\hline
 $b^{(0)}_0$    &$2-2a_2$& $-2-\frac{3a_4}{4}$&$-\frac16-\frac{5a_4}{3}$ &$\frac{13}2-3a_4$& $-2-\frac{19a_5}{18}$ &$\frac76-\frac{19a_5}9$ \\\hline\hline
 $b_{\frac12}$  & $a_2 $ & $a_4$              & $a_4 $                   & $ a_4 $         & $a_5$ &$a_5$             \\\hline\hline
 $b^{(0)}_1$    &        & $2-\frac{3a_4}{4}$ &   $ \frac16 (5-a_4) $    &                 & $2-\frac{5a_5}{9}$ &$\frac{6-a_5}{9}$ \\\hline
 $b_{\frac32}$  &        &$\frac{-2+a_4}{4}$  &                          &                 & $\frac{a_5-3}{6}$& \\\hline\hline
 %$b^{(0)}_2$    &        &                    &                          &                 &\\\hline
 %$b_{\frac52}$  &         &                   &                          &                 &\\\hline\hline
 CFL$_\text{max}$ (RK3)&1& 0.7985             &    1                     & 0.45            & 0.675 &0.855 \\
 parameter      &$a_2\in[1,2]$&$a_4\simeq 1.7723$& $a_4\in[\frac12,1.1]$ & $a_4\simeq 3.5$ & $a_5 \simeq 1.6$ &$a_5 \simeq 1.5$\\
 region & Fig. \ref{fig:rk3stabilityregions1st} & Fig. \ref{fig:rk3stabilityregions3rd} & Fig. \ref{fig:rk3stabilityregions3upwind1}& Fig. \ref{fig:rk3stabilityregions3upwind2}& Fig. \ref{fig:rk3stabilityregions4} & Fig. \ref{fig:rk3stabilityregions4upwind1}
\end{tabular}
\caption{Overview of finite difference formulae. The lower part shows the maximum attainable CFL number together with the corresponding value of the parameter. The formula FD2 becomes FD3 upon the choice $a_2=4$. }
\label{tab:finitedifferences}
\end{table}

\begin{table}
 \centering
\begin{tabular}{c||c|c|c||c||c|c|}
 name           &  FD6a             & FD6b                 & FD6c             & FD7 & FD8a                          & FD8c  \\
 $p$            &  6                &        6             &   6              &7&          8                        & 8\\
 \hline\hline
 $b^{(0)}_{-3}$ &                   &                      &                  & &                                   &\\\hline
 $b_{-\frac52}$ &                   &                      &                  & & $\frac{-8+3 a_8}{48}$             &\\\hline
 $b^{(0)}_{-2}$ &                   &                      &$\frac{4-a_6}{12}$& $\frac{2-a_7}{48}$ & $\frac{49}{72} -\frac{25a_8}{96}$ &\\\hline
 $b_{-\frac32}$ &                   &  $\frac{a_6-1}{9}$   &$ a_6-\frac{11}3$ & $\frac{3a_7-5}{9}$& $ -\frac{7}{3} + a_8 $            &$ \frac{1+a_8}{36}$\\\hline
 $b^{(0)}_{-1}$ &$-\frac{1+a_6}{36}$&$\frac{19-22a_6}{54}$ &$\frac{302-87a_6}{36}$& $\frac{586-393a_7}{432}$ &$\frac{293}{72}-\frac{185a_8}{96}$ &$ -\frac{28+25a_8}{216}$\\\hline
 $b_{-\frac12}$ &$\frac{2+a_6}{3}$  &  $a_6$               &$ 3a_6-8$         & $\frac{3a_7}2-1$& $ -3 +\frac{9a_8}{4} $            &$ \frac49 (2+a_8)$\\\hline
 $b^{(0)}_0$    &$-\frac94-\frac{29a_6}{36}$&$-\frac{89+76a_6}{54}$& $ \frac{86-87a_6}{36}$ & $\frac{-494-717a_7}{432}$ &$-\frac{31}{72}-\frac{185a_8}{96}$ & $-\frac{5}{216}(108+37a_8)$\\\hline\hline
 $b_{\frac12}$  &     $a_6$         & $ a_6 $              &$a_6$             & $a_7$ & $ a_8 $                           &$a_8$\\\hline\hline
 $b^{(0)}_1$    &$\frac94-\frac{29a_6}{36}$&$\frac{50 -11a_6}{27}$&$\frac59-\frac{a_6}{12}$& $\frac{730-141a_7}{432}$& $\frac{436-75a_8}{288} $          &$ \frac52 - \frac{185a_8}{216}$\\\hline
 $b_{\frac32}$  & $\frac{a_6-2}{3}$ &  $ \frac{a_6-4}{9} $ &                  & $\frac{3a_7-14}{36}$ &$-\frac13+\frac{a_8}{16}$           &$ \frac49 (a_8-2)$\\\hline
 $b^{(0)}_2$    & $\frac{1-a_6}{36}$&                      &                  & &                                    &$ \frac{28-25 a_8}{216}$\\\hline
 $b_{\frac52}$  &                   &                      &                  & &                                   & $\frac{a_8-1}{36}$\\\hline\hline
 CFL$_\text{max}$ (RK3)& $0.7$      &  $0.713$             &$0.56$            &$0.73 $& $0.657$                     & $0.62$ \\
 parameter      &$a_6\simeq1.88$    &$a_6\simeq\frac14$    &$a_6 \simeq 2.3$  &$a_7\simeq 0.68$& $a_8 \simeq \frac43$        &  $a_8\simeq 1.9$ \\
 region & Fig. \ref{fig:rk3stabilityregions5thupwindminus1}& Fig. \ref{fig:rk3stabilityregions5th} & & Fig. \ref{fig:rk3stabilityregions6th} & Fig. \ref{fig:rk3stabilityregions7th} & Fig. \ref{fig:rk3stabilityregions7thupwindminus1}
\end{tabular}
\caption{Table \ref{tab:finitedifferences} continued.}
\label{tab:finitedifferencescont}
\end{table}

\subsection{Transonic upwinding}\label{ssec:findifftransonicupwinding}

  By analogy with the way how transonic upwinding has been proposed to be handled in \cite{barsukow19activeflux} (see also Section \ref{ssec:transonicupwindinggeneral}), we propose here for the case of a scalar conservation law ($m=1$) the following choice:

\begin{align}
 \tilde q_{i+\frac12} = \begin{cases}
 q^n_{i+\frac12} & \text{if } \sign f'(q^n_{i-\frac12}) = \sign f'(q^n_{i+\frac12}) = \sign f'(q^n_{i+\frac32})\\
 q^n_{i+\frac12} & \text{else and if } q^n_{i-\frac12} < q^n_{i+\frac32} \text{transonic rarefaction}\\
 q^n_{i-\frac12} & \text{else and if } |f'(q^n_{i-\frac12})| \geq \max(|f'(q^n_{i+\frac12})|, |f'(q^n_{i+\frac32})|) \\
 q^n_{i+\frac12} & \text{else and if } |f'(q^n_{i+\frac12})| \geq \max(|f'(q^n_{i-\frac12})|, |f'(q^n_{i+\frac32})|) \\
 q^n_{i+\frac32} & \text{else and if } |f'(q^n_{i+\frac32})| \geq \max(|f'(q^n_{i-\frac12})|, |f'(q^n_{i+\frac12})|) \\
                        \end{cases}
\end{align}

While the original modification from \cite{barsukow19activeflux} leads to a third-order accurate method, the above suggestion degrades the accuracy at sonic points. We do not find averages such as $\tilde q_{i+\frac12} = \frac{q_{i-\frac12}^n + 2 q_{i+\frac12}^n + q_{i+\frac32}^n}{4}$ to work satisfactorily for our high order versions.

For systems, the presence of different waves and thus inherent information transport from different parts of the initial data seems to make transonic upwinding less important. We experimentally find that the high-order version of Active Flux via finite differences in many cases, such as in the numerical examples below, works well with the simple choice $\tilde q_{i+\frac12} = q^n_{i+\frac12}$.

\subsection{Limiting}

We propose to gradually reduce the order of the finite difference approximation employed if a violation of monotonicity is detected. To this end, at every cell interface, the candidate derivative approximations $D$ and $D^*$ are first computed using the finite differences of the desired order of accuracy of the method. Then, for each of these formulae the following algorithm is performed:

\begin{figure}
 \centering
 \includegraphics[width=0.35\textwidth]{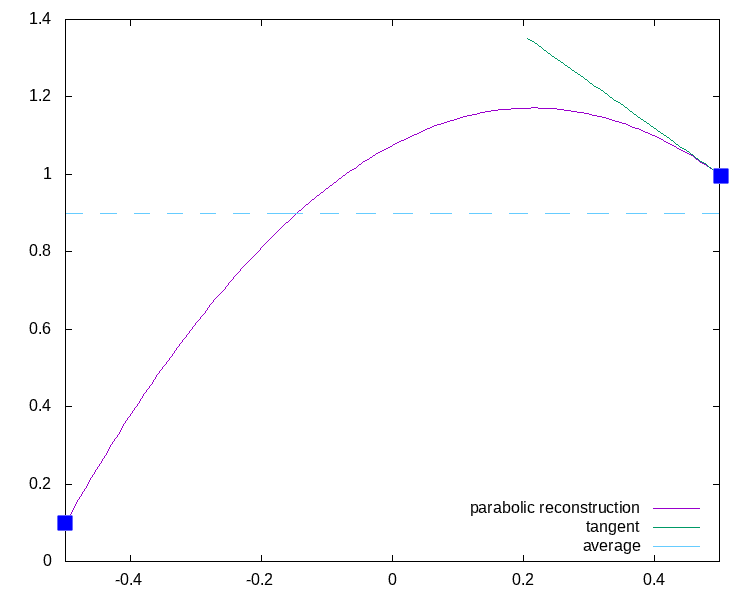}
 \caption{The three values $q_{i-\frac12} < \bar q_i < q_{i+\frac12}$ are monotone and growing, but the finite difference formula obtained from differentiating the reconstruction at the right edge of the interval obtains a negative slope. This is a situation that the limiting procedure is correcting.}
 \label{fig:findiffnonmonotonemodif}
\end{figure}

\begin{enumerate}[1.]
\item We check whether the values
\begin{align}
 \bar q_{i-\ell}, q_{i+\frac12-\ell}, \ldots, \bar q_{i+\ell}, q_{i+\frac12+\ell} \label{eq:valuesinstencil}
\end{align}
involved in the finite difference formula are monotone. If they are not, no further action is taken and the finite difference approxmation is accepted.
\item If the values in the finite difference are monotone, then we check whether the finite difference approximation to the derivative has the same sign as, say, $q_{i+\frac12} - \bar q_i$ (compare to Figure \ref{fig:findiffnonmonotonemodif}). If so, the finite difference approximation is accepted.
\item If the finite difference has the opposite sign, we reduce the order by one and repeat steps 1 and 2, descending e.g. as FD8a--FD7--FD6b--FD5b--FD4b--FD3 (that is, FD2 with $a_2=4$), until either monotonicity is obtained, or until even the finite difference of lowest order (second) is found to be non-monotone. 
\item If all the finite differences including FD3 are found to be non-monotone, a fall-back finite difference is used. We propose to take inspiration from the power-law limiting suggested in \cite{barsukow19activeflux}. Its original version amounts to a modification of the reconstruction, but a finite difference can be obtained by differentiation at the cell interfaces $x = \pm \frac{\Delta x}{2}$. The nonvanishing derivatives of the two power laws \eqref{eq:limitingpowerlaw1}--\eqref{eq:limitingpowerlaw2} yield two nonlinear finite difference approximations:
\begin{align}
 \frac{\dd}{\dd x} q^n_{\text{recon},i,\text{power-law},1}\left(-\frac{\Delta x}{2}\right) &= \frac{1}{\Delta x} (q_{i+\frac12}^n - q_{i-\frac12}^n) \frac{\bar q_i^n - q_{i-\frac12}^n}{q_{i+\frac12}^n - \bar q_i^n} \label{eq:powerlawderivative1}\\
 \frac{\dd}{\dd x} q^n_{\text{recon},i,\text{power-law},2}\left(\frac{\Delta x}{2}\right) &= \frac{1}{\Delta x}(q_{i+\frac12}^n - q_{i-\frac12}^n) \frac{q_{i+\frac12}^n - \bar q_i^n}{\bar q_i^n - q_{i-\frac12}^n} \label{eq:powerlawderivative2}
\end{align}
These finite differences are first-order accurate, as can be checked explicitly. Up to a shift by one cell, the one is a finite difference at $x_{i+\frac12}$ for positive speed, and the other for negative, i.e. we define
\begin{align}
 D_1(q_{i-\frac12}, \bar q_i, q_{i+\frac12}) &:= \frac{1}{\Delta x}(q_{i+\frac12}^n - q_{i-\frac12}^n) \frac{q_{i+\frac12}^n - \bar q_i^n}{\bar q_i^n - q_{i-\frac12}^n} \\
 D_1^*(q_{i+\frac12}, \bar q_{i+1}, q_{i+\frac32}) &:= \frac{1}{\Delta x} (q_{i+\frac32}^n - q_{i+\frac12}^n) \frac{\bar q_{i+1}^n - q_{i+\frac12}^n}{q_{i+\frac32}^n - \bar q_{i+1}^n}  
\end{align}
\end{enumerate}
By analogy with the original version from \cite{barsukow19activeflux}, and as said above, this finite difference is used only if the values \eqref{eq:valuesinstencil} are monotone and if all higher order finite differences were found to violate monotonicity. As the second order finite difference is among them, this in particular means that \eqref{eq:powerlawderivative1}--\eqref{eq:powerlawderivative2} are only used if condition \eqref{eq:monotnicityconditionparabola} is violated. This implies that they are well-behaved. For practical reasons we also refrain from using \eqref{eq:powerlawderivative1}--\eqref{eq:powerlawderivative2} if the exponents of the power laws are outside the interval $[\frac{1}{50},50]$, and use the finite difference FD3 instead.

%As has been pointed out in \cite{abgrall20}, it is possible to use different variables for the cell averages and for the point values. For the Euler equations, for example, one could use conserved variables (density, momentum, energy) for the cell averages, and primitive variables for the point values. The latter are much easier to handle, and are often of greater practical interest in simulations. In the examples of systems below, however, all degrees of freedom employ the conservative variables, and the same is true for the reconstruction.

\subsection{Numerical results}

\subsubsection{Linear advection}

\begin{figure} % cb18c09
  \centering
  \includegraphics[width=0.49\textwidth]{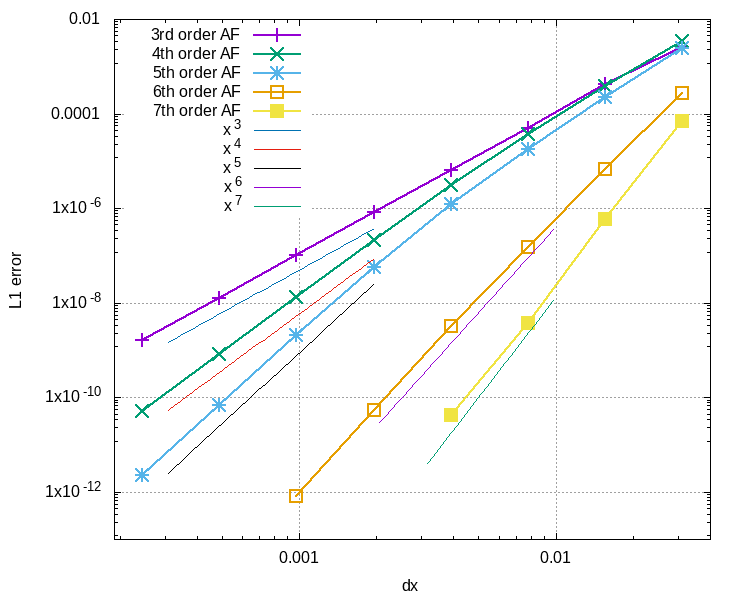}\hfill\includegraphics[width=0.49\textwidth]{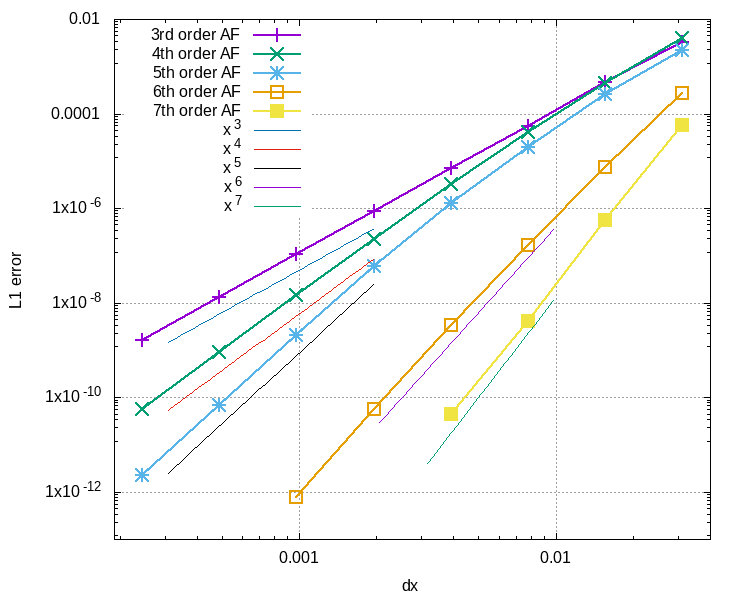}
  \caption{Convergence of the the high order extension of Active Flux via finite differences for linear advection. \textit{Left}: $L^1$ error of the point values. \textit{Right}: $L^1$ error of the averages.}
  \label{fig:findiff-advection-convergence}
 \end{figure}

Figure \ref{fig:findiff-advection-convergence} shows the experimental order of accuracy for high order versions of Active Flux using finite differences for linear advection $\del_t q + c \del_x q= 0$ ($c = 1$). Numerical methods of orders 3 to 7 are shown. The tests are using a Runge-Kutta integrator of 3rd order accuracy and a small CFL number of $10^{-2}$ in order for the total error to be dominated by the spatial error. The initial data are
\begin{align}
 q_0(x) = 0.8 + \exp\left( - \frac{(x - 0.5)^2}{0.05^2} \right) \label{eq:initialconvergence}
\end{align}
on grids covering $[0,1]$. The error is shown at time $t=0.1$. Limiting is not used. The finite differences employed are 
FD7 ($a_7 =2.5$), FD6b ($a_6 = 2$), FD5b ($a_5 = 1.55$), FD4b ($a_4 = 1$), FD3 (FD2 with $a_2=4$).
We find that it is not always easy in practice to demonstrate experimentally the desired order of accuracy if the parameters of the finite difference formulas are chosen too close to the stability limit. We also find that the choice of the stencil among those of the same order of accuracy can affect the value of the error; a detailed study of this influence shall, however, remain subject of future work. Generally, we find the theoretical orders of accuracy reflected in the simulation results, as Figure \ref{fig:findiff-advection-convergence} shows.

\subsubsection{Burgers' equation} \label{ssec:findiffburgers}

\begin{figure} % e72c08e
  \centering
\includegraphics[width=0.49\textwidth]{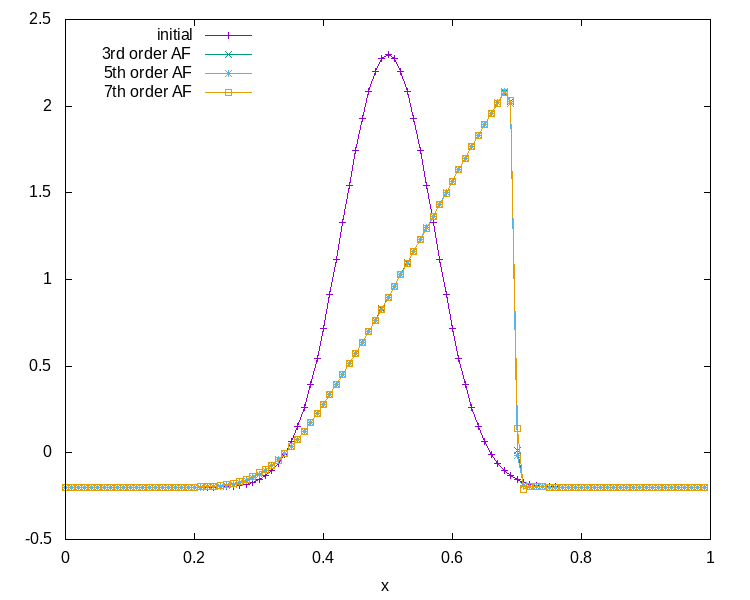}\hfill\includegraphics[width=0.49\textwidth]{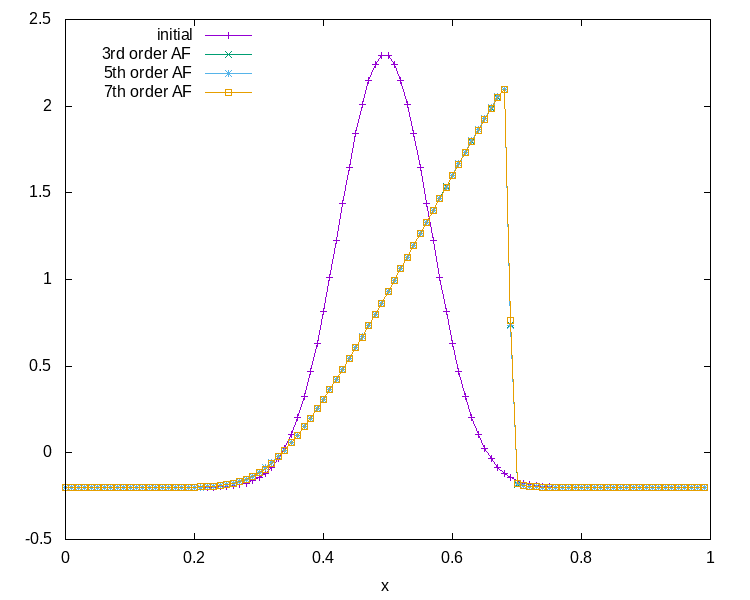}
  \caption{Numerical evolution of Gaussian initial data for Burgers' equation at time $t=0.1$. \textit{Left}: Point values. \textit{Right}: Averages.}
  \label{fig:findiff-burgers}
 \end{figure}

Figure \ref{fig:findiff-burgers} shows the numerical time evolution of Gaussian initial data
\begin{align}
 q_0(x) = 2.5 \exp\left( - \frac{(x - 0.5)^2}{0.1^2} \right) - 0.2 \label{eq:initialtransonicburgers}
\end{align}
on a grid of 100 cells covering $[0,1]$ subject to Burgers' equation. This setup is transonic, but is captured without artifacts due to use of limiting and transonic upwinding as described above. The CFL number used is 0.4.
 
\begin{figure} % 495a672
  \centering
  \includegraphics[width=0.49\textwidth]{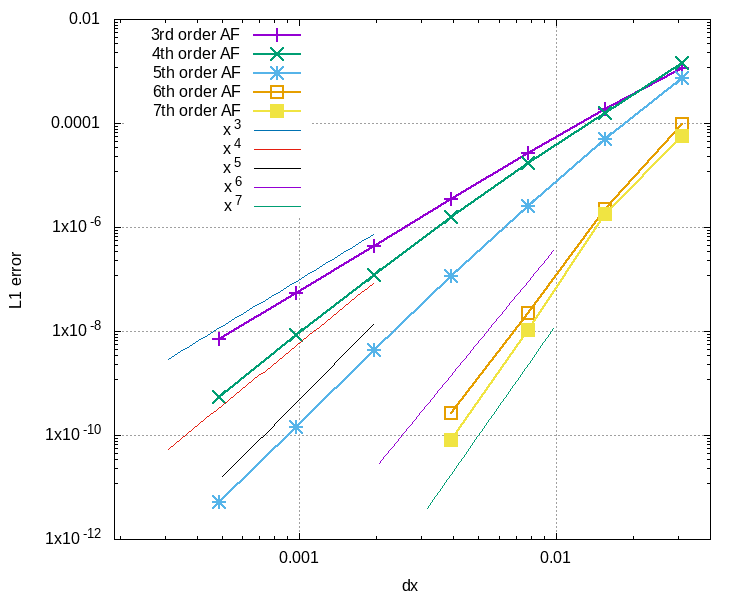}\hfill\includegraphics[width=0.49\textwidth]{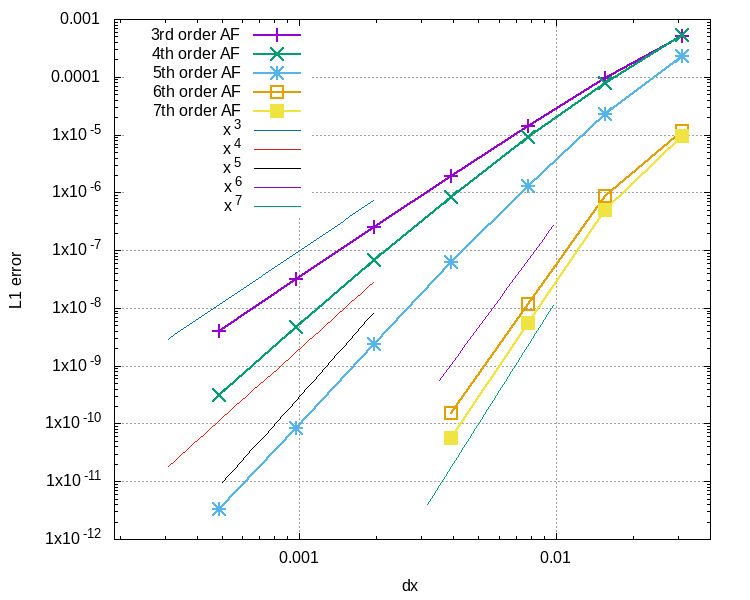}
  \caption{Convergence of the the high order extension of Active Flux via finite differences for Burgers' equation. \textit{Left}: $L^1$ error of the point values. \textit{Right}: $L^1$ error of the averages.}
  \label{fig:findiff-burgers-convergence}
 \end{figure}

Figure \ref{fig:findiff-burgers-convergence} shows the experimental convergence on different grids covering $[0,1]$ using initial datum \eqref{eq:initialconvergence} at time $t=0.01$, before the formation of a shock. The CFL number is chosen to be $10^{-3}$, to ensure that the error is dominated by that of spatial discretization. Here, limiting is not used. 
 
\subsubsection{Euler equations}

Finally, we assess the performance of the method on the Euler equations
\begin{align}
 \del_t \rho + \del_x (\rho v) &= 0 & \rho & \colon \mathbb R^+_0 \times \mathbb R \to \mathbb R^+\\
 \del_t (\rho v) + \del_x (\rho v^2 + p) &= 0 & v & \colon \mathbb R^+_0 \times \mathbb R \to \mathbb R\\
 \del_t e + \del_x (\rho v) &= 0& e & \colon \mathbb R^+_0 \times \mathbb R \to \mathbb R^+
\end{align}
with
\begin{align}
 e &= \frac{p}{\gamma-1} + \frac12 \rho v^2
\end{align}

\begin{figure} % 45071c2
  \centering
\includegraphics[width=0.49\textwidth]{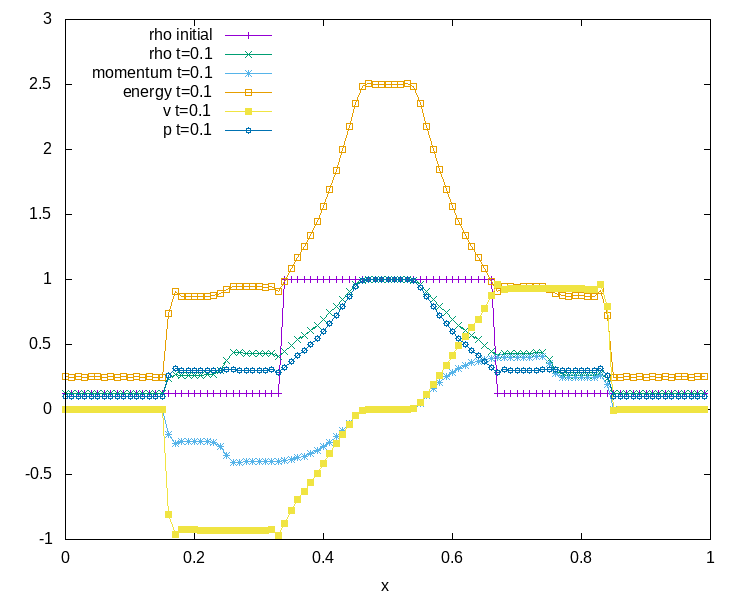}\hfill\includegraphics[width=0.49\textwidth]{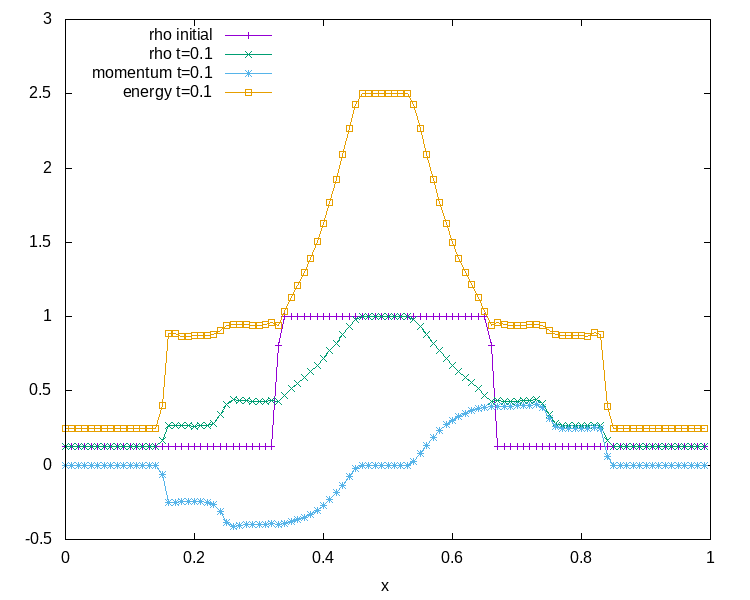}
  \caption{Numerical evolution of the Sod shock tube at time $t=0.1$. \textit{Left}: Point values. \textit{Right}: Averages.}
  \label{fig:findiffeulersod7}
 \end{figure}

Figure \ref{fig:findiffeulersod7} shows the results for the Sod shock tube \cite{sod78} on a grid of $100$ points; the CFL number is 0.25 and the 6th order method is used.

\section{High order via further point values} \label{sec:points}

A different viewpoint is that of Active Flux as an enriched Finite Volume method, which includes point values additionally to the cell average. Higher order can thus be obtained by including more point values, an idea also mentioned in \cite{he21}. In one dimension, the cell interfaces are already ``taken'' and any new point values must be located inside the cell. With $k$ such points, the generalization of the reconstruction \eqref{eq:parabolicrecon} in cell $i$ then is a polynomial $q_{\text{recon},i} \in P^{k+2}$, while the global regularity remains $C^0$. Therefore, the order of accuracy is expected to be $\order = k+3$, while the average number of degrees of freedom per cell is $k+2$. Denote the locations of the additional points by $\Delta x \xi_j$, $\xi_j \in \left(-\frac12,\frac12\right)$, and their values by $q^n_{i,j}$, $j = 1,\ldots,k$.

Such a method is thus a natural extension of the traditional Active Flux method (Definition \ref{def:activefluxtraditional}) and it is natural to evolve the new points by using an evolution operator. First, all point values are updated, and then the average is updated with fluxes that are obtained through quadrature in time at the cell interfaces. In order to achieve a high order quadrature, point values at cell interfaces are computed at intermediate times as well, while the point values inside the cell are evolved over the full time step directly. The evolution operators, as well as the limiting strategy and the transonic upwinding as described in \cite{barsukow19activeflux} and summarized in Sections \ref{ssec:afwithevolutionop}, \ref{ssec:limiting}--\ref{ssec:transonicupwindinggeneral} require only the smallest modifications.

\subsection{High-order reconstruction}

The place of the traditional reconstruction \eqref{eq:recon3rdstandard} is taken by a reconstruction in cell $i$
\begin{align}
q_{\text{recon},i} &\in C^0 \cap L^1_\text{loc} & q_{\text{recon},i} &\colon (\mathbb R^m)^{3+k} \times \left[-\frac{\Delta x}{2},\frac{\Delta x}{2}\right] \to \mathbb R^m \label{eq:reconpointsdef}
\end{align}
satisfying
\begin{align}
q_{\text{recon},i}\left(q_{i-\frac12}, \bar q_i, q_{i+\frac12}, \pm \frac{\Delta x}{2}\right ) &= q_{i\pm\frac12}\\
q_{\text{recon},i}\left(q_{i-\frac12}, \bar q_i, q_{i+\frac12},  \Delta x \xi_j \right ) &= q_{i,j} \qquad j =1, \ldots, k\\
 \frac{1}{\Delta x} \int_{-\frac{\Delta x}{2}}^{\frac{\Delta x}{2}} \dd x \, q_{\text{recon},i}\left(q_{i-\frac12}, \bar q_i, q_{i+\frac12}, x\right ) &= \bar q_i \label{eq:reconpointsaverage}
\end{align}

%\begin{align}
% q_\text{recon}&\left(\phi(x_{i-\frac12}), \frac{1}{\Delta x} \int_{x_{i-\frac12}}^{x_{i+\frac12}} \dd x \,\phi(x), \phi(x_{i+\frac12}), \phi(x_{i} + \Delta x \xi_1), \ldots, \phi(x_i + \Delta x \xi_k),   x\right ) \label{eq:reconpoints}\\&\phantom{mmmmmmmmmmmmmmmmm}\nonumber= \phi(x_i + x) + \mathcal O(\Delta x^{3+k}) \qquad \forall x \in \left[-\frac{\Delta x}{2},\frac{\Delta x}{2}\right]
%\end{align}
%for any analytic function $\phi \colon \mathbb R \to \mathbb  R$. 

These are $k+3$ conditions. It is natural to seek the reconstruction in the space $P^{2+k}$ of polynomials which will yield an Active Flux method of $\order$-th order with $k+2 = \order-1$.
%\begin{align}
% q^n_{\text{recon},i}\left( \pm \frac{\Delta x}{2}  \right ) &= q^n_{i\pm\frac12} & 
% \frac{1}{\Delta x} \int_{-\frac{\Delta x}{2}}^{\frac{\Delta x}{2}} q^n_{\text{recon},i}(x) \,\dd x &= \bar q_i^n \label{eq:reconpoints1} \\
% q^n_{\text{recon},i}\left( \Delta x \xi_j \right ) &= q^n_{i,j} \qquad j = 1, \ldots, k \label{eq:reconpoints2}
%\end{align}
In this case, it is useful to rewrite this interpolation problem (which involves interpolating point values and averages) as a standard pointwise interpolation problem, in order to be able to use standard algorithms.

\begin{theorem} \label{thm:highorderinterpolation}
 Let $\{ (x_j,y_j) \in \mathbb R^2 \}_{j=1, \ldots, \order-1}$ with $(x_i = x_j) \Rightarrow (i = j) \, \forall i,j$ and let $\bar y \in \mathbb R$ be given. Fix an arbitrary $x_0 \neq x_j \,\forall j$. Then there exist unique polynomials
 \begin{align}
 p_1 &\in P^{\order-2} \text{ with } p(x_j) = y_j \quad \forall j = 1, \ldots, \order-1 \label{eq:thmpoints1}\\
 p_2 &\in P^{\order-1} \text{ with } \begin{cases} p(x_j) = 0 & \forall j = 1, \ldots, \order-1\\ p(x_0) = 1 \end{cases} \label{eq:thmpoints2}
\end{align}  
 Moreover, if it exists, then the unique polynomial $p \in P^{k+2}$ fulfilling 
 \begin{align}
  p(x_j) &= y_j \qquad j = 1, \ldots \order-1 & &\text{and} & \frac{1}{\Delta x} \int_{-\frac{\Delta x}{2}}^{\frac{\Delta x}{2}} \dd x \, p(x) &= \bar y
 \end{align}
 is 
 \begin{align}
 p(x) = p_1(x) + p_2(x)  \cdot \frac{\Delta x \bar y - \int_{-\frac{\Delta x}{2}}^{\frac{\Delta x}{2}}\dd x \, p_1 }{\int_{-\frac{\Delta x}{2}}^{\frac{\Delta x}{2}} \dd x \, p_2}
\end{align}
\end{theorem}

\begin{proof}
$p_1$ and $p_2$ are defined as standard interpolating polynomials, they thus exist and are unique. By construction, for any $\alpha \in \mathbb R$
\begin{align}
 p(x_j) = p_1(x_j) + \alpha p_2(x_j) = y_j \quad \forall j = 1, \ldots, \order-1
\end{align}
and upon integration
\begin{align}
 \Delta x \bar y = \int_{-\frac{\Delta x}{2}}^{\frac{\Delta x}{2}} \dd x \, p_\text{recon} = \int_{-\frac{\Delta x}{2}}^{\frac{\Delta x}{2}}\dd x \, p_1  + \alpha \int_{-\frac{\Delta x}{2}}^{\frac{\Delta x}{2}} \dd x \, p_2
\end{align}
i.e.
\begin{align}
 \alpha = \frac{\Delta x \bar y - \int_{-\frac{\Delta x}{2}}^{\frac{\Delta x}{2}}\dd x \, p_1}{\int_{-\frac{\Delta x}{2}}^{\frac{\Delta x}{2}} \dd x \, p_2}
\end{align}
if $\int_{-\frac{\Delta x}{2}}^{\frac{\Delta x}{2}} \dd x \, p_2 \neq 0$.
\end{proof}

This theorem can be used to define the high-order reconstruction through following identifications:
\begin{align}
 q^n_{\text{recon},i}(x) := p(x)
\end{align}
with 
\begin{align}
x_1 &:= -\frac{\Delta x}{2}   &  x_j &:= \Delta x \xi_{j-1} \quad j = 2, \ldots, \order-2 & x_{\order-1} &:= \frac{\Delta x}{2}\\
y_1 &:= q^n_{i-\frac12} &    y_j &:= q^n_{i,j-1} \quad j = 2,\ldots,\order-2 & y_{\order-1} &:= q_{i+\frac12}^n
\end{align}

\begin{theorem}
 Consider the setup of Theorem \ref{thm:highorderinterpolation} and assume the points $\{x_j\}_{j=1, \ldots, \order-1} $ to be located symmetrically around 0. If their number $\order-1$ is odd, then $\int_{-\frac{\Delta x}{2}}^{\frac{\Delta x}{2}} \dd x \, p_2 = 0$, and there is no polynomial of degree $k+2$ fulfilling \eqref{eq:thmpoints1}--\eqref{eq:thmpoints2}.
\end{theorem}
\begin{proof}
Observe that, with the Lagrange formula, $p_2$ can be given the explicit form
\begin{align}
 p_2(x) = \prod_{j=1}^{\order-1} \frac{x - x_j}{x_0 - x_j}
\end{align}
Take $n$ to be even, and call $\nu := \order/2$. Order the points symmetrically as 
\begin{align}
x_1 < x_2 < \ldots < x_{\nu-1} < x_{\nu} = 0 < -x_{\nu-1} < \ldots < -x_2 < -x_{1}
\end{align}

Then, collecting the pairs,
\begin{align}
 p_2(x) = \frac{x}{x_0 - x_\nu}  \prod_{j=1}^{\nu-1} \frac{x^2 - x_j^2}{x_0^2 - x_j^2} 
\end{align}
Thus, $p_2$ is antisymmetric, and $\int_{-\frac{\Delta x}{2}}^{\frac{\Delta x}{2}} \dd x \, p_2 = 0$.
\end{proof}
This means that we are unable to construct methods of even order with the point values located symmetrically around the cell center. 
There exist other numerical methods where odd orders of accuracy are given preference, e.g. WENO. If one would be interested to construct an even-order version of Active Flux, one can show that this is possible upon usage of an asymmetric distribution of points. 

%This is simiar to Roe's choice of storing derivatives at cell interfaces (Section \ref{}), where the method naturally has an odd order as well.

\begin{definition}[Active Flux with additional point values]
With the reconstruction $q_{\text{recon},i}$ in cell $i$ defined in \eqref{eq:reconpointsdef}--\eqref{eq:reconpointsaverage}, define the global reconstruction 
\begin{align}
q_\text{recon} &\in  C^0 \cap L^1_\text{loc} & q_\text{recon} &\colon \mathbb R \to \mathbb R^m 
\end{align}
\begin{align} 
 q_\text{recon}(x) &:= q_{\text{recon},i}\left(q_{i-\frac12}, \bar q_i, q_{i+\frac12}, q_{i,1}, \ldots, q_{i,k}, x-x_i\right ) \quad \text{ if }x\in [x_{i-\frac12},x_{i+\frac12}] \label{eq:pointsglobalrecon}
\end{align}

The \emph{Active Flux method with additional point values} is the following semi-discretization
\begin{align}
 \begin{cases} \displaystyle \frac{\dd}{\dd t} \bar q_i(t) = - \frac{f(q_{i+\frac12}(t)) - f(q_{i-\frac12}(t))}{\Delta x} \\ \\
 \displaystyle q_{i+\frac12}(t) = \Big ( \text{solution at $x=x_{i+\frac12}$ of the IVP \eqref{eq:conslaw} with initial data } q_{\text{recon}} \Big ) + \mathcal O(t^{3+k}) \\ \\
 \displaystyle q_{i,j}(t) = \Big ( \text{solution at $x=x_{i} + \Delta x \xi_j$ of the IVP \eqref{eq:conslaw} with initial data } q_{\text{recon}} \Big ) + \mathcal O(t^{3+k}) \\ \nonumber \phantom{mmmmmmmmmmmmmmmmmmmmmmmmmmmmmmm} j = 1,\ldots,k
 \end{cases} \label{eq:activefluxdefpoints}
\end{align}
of \eqref{eq:conslaw} with the interpretations
\begin{align}
 \bar q_i(t) &\simeq \frac{1}{\Delta x} \int_{x_{i-\frac12}}^{x_{i+\frac12}} \dd x \,  q(t, x) &
 q_{i+\frac12}(t) &\simeq q(t, x_{i+\frac12})\\
 && q_{i,j}(t) &\simeq q(t, x_{i} + \Delta x \xi_j) \qquad j = 1, \ldots, k
\end{align}
\end{definition}

A solution operator for linear advection $\del_t q + c \del_x q = 0$ is easily given:
\begin{align}
 q_{i+\frac12}(t) = q_\text{recon}(x_{i+\frac12} - c t) &= \begin{cases}  q_{\text{recon},i}(\Delta x/2 - c t) & c > 0, ct < \Delta x \\
                                                            q_{\text{recon},i+1}(-\Delta x/2 - c t) & c < 0, |c|t < \Delta x 
                                                           \end{cases}
\end{align}
%The choice of the cell shall be encoded by defining a function $\sigma \in \{0, 1\}$, such that one writes (still assuming $|c|t < \Delta x$)
%\begin{align}
% q_{i+\frac12}(t) = q_{\text{recon},i + \sigma}(\Delta x/2 - \sigma \Delta x - c t) \qquad \sigma = \begin{cases} 0 & c > 0 \\ 1 & \text{else} \end{cases}
%\end{align}

A solution operator for scalar conservation laws based on a fixpoint iteration to find the characteristic is given in \cite{barsukow19activeflux} for any order of accuracy. To this end, the equation for the foot point $x_0$ of the characteristic
\begin{align}
  x_0 + f'(q_\text{recon}(x_0)) t = x_{i+\frac12}
\end{align}
is solved iteratively:
\begin{align}
 x_0^{(0)} &:= x_{i+\frac12}  \label{eq:fixpointiterationinit}\\ 
 x_0^{(\ell)} &:= x_{i+\frac12} - f'(q_\text{recon}(x_0^{(\ell-1)})) t \qquad \ell = 1, \ldots, \ell_\text{max} \label{eq:fixpointiteration}
\end{align}
An approximate solution operator is
\begin{align}
 q_{i+\frac12}(t) \simeq q_\text{recon}(x_0^{(\ell_\text{max})}) = q(t, x_{i+\frac12}) +  \mathcal O(t^{\ell_\text{max}+1}) \label{eq:fixpointiterationerror}
\end{align}

A third order accurate evolution operator for hyperbolic systems of conservation laws can also be found in \cite{barsukow19activeflux}, and an approximate evolution operator for inhomogeneous problems has been given in \cite{barsukow19activefluxsource},\cite{barsukow20swaf}.

\subsection{Time integration and stability conditions}

As in Section \ref{ssec:afwithevolutionop}, the high-order extension via additional point values can be immediately formulated as an explicit method. For linear advection $\del_t q + c \del_x q = 0$ with $c> 0$, the algorithm extending \eqref{eq:updatelinearadvectionvanleer}--\eqref{eq:updateaverageleapfrogvanleer} thus reads

%\begin{align}
% q_{i+\frac12}^{n+\frac{\ell}{M-1}} &= q^n_{\text{recon},i}(x_{i+\frac12} - c \zeta_\ell \Delta t) \qquad \ell = 1,\ldots,M-1\\
% q_{i,j}^{n+1} &=  q^n_{\text{recon},i - \sigma_j(i)}(\Delta x (\sigma_j(i) +\xi_j) - c \Delta t) = 0 \qquad j = 1,\ldots,\order-3\\
% \bar q_i^{n+1} &= \bar q_i^n - \Delta t \sum_{\ell = 0}^{M-1} \omega_\ell \frac{f(q^{n + \frac{\ell}{M-1}}_{i+\frac12}) - f(q^{n + \frac{\ell}{M-1}}_{i-\frac12})}{\Delta x}   
%\end{align}

\begin{align}
 q_{i+\frac12}^{n+\frac{\ell}{M-1}} &= q^n_\text{recon}(x_{i+\frac12} - c \zeta_\ell \Delta t) \qquad \ell = 1,\ldots,M-1\\
 q_{i,j}^{n+1} &=  q^n_\text{recon}(x_{i} + \Delta x \xi_j - c \Delta t)  \qquad j = 1,\ldots,\order-3 \label{eq:pointsupdateintermedaite}\\
 \bar q_i^{n+1} &= \bar q_i^n - \Delta t \sum_{\ell = 0}^{M-1} \omega_\ell \frac{f\left(q^{n + \frac{\ell}{M-1}}_{i+\frac12}\right) - f\left(q^{n + \frac{\ell}{M-1}}_{i-\frac12}\right)}{\Delta x}   
\end{align}

Here, $M$ is the number of points in the time-quadrature approximating the time-step-averaged flux
\begin{align}
\frac{1}{\Delta t} \int_0^{\Delta t} f(q_{i+\frac12}(t)) \dd t \simeq \sum_{\ell = 0}^{M-1} \omega_\ell f(q_{i+\frac12}^{n+\frac{\ell}{M-1}}) \label{eq:pointsfluxquadrature}
\end{align}
(with $\zeta_0 = 0$, $\zeta_{M-1} = 1$), and $\{\omega_\ell\}_{\ell = 0, \ldots, M-1}$ the corresponding weights. Clearly, $q_{i+\frac12}^{n+\frac{\ell}{M-1}} \simeq q_{i+\frac12}(\Delta t \zeta_\ell)$. Note that we use the notation $q_{i+\frac12}^{n+\frac{\ell}{M-1}}$ to denote intermediate values, but we do not actually insist on them being equidistant in time. In fact, it is natural to use a Gauss-Lobatto quadrature here. 

\begin{figure}
 \centering
 \includegraphics[width=0.32\textwidth]{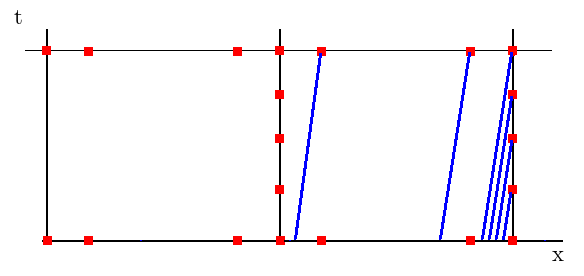} \hfill \includegraphics[width=0.32\textwidth]{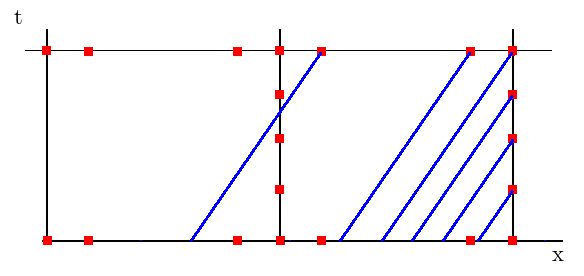}  \hfill \includegraphics[width=0.32\textwidth]{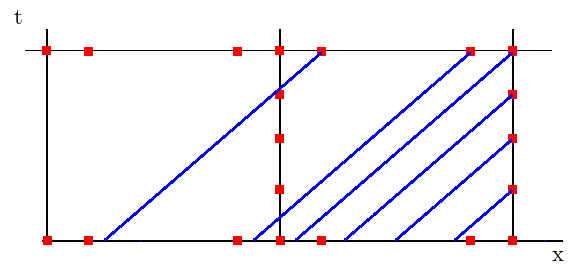}
 \includegraphics[width=0.32\textwidth]{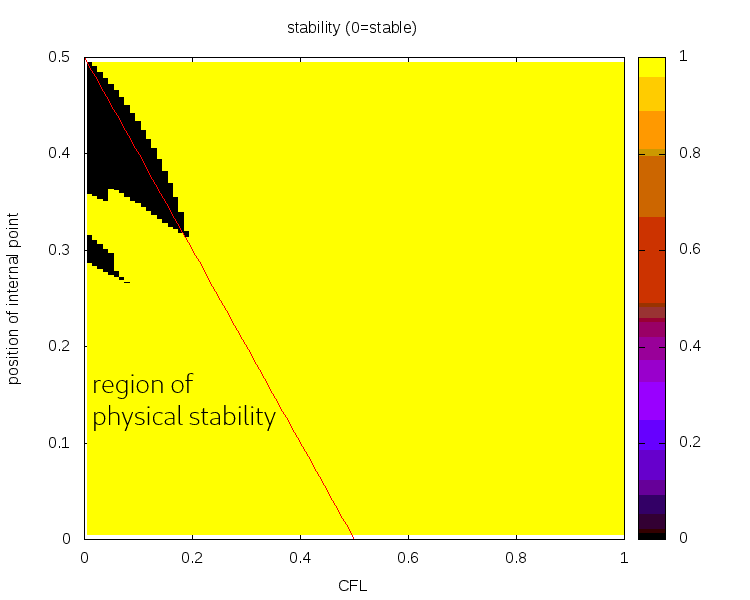} \hfill \includegraphics[width=0.32\textwidth]{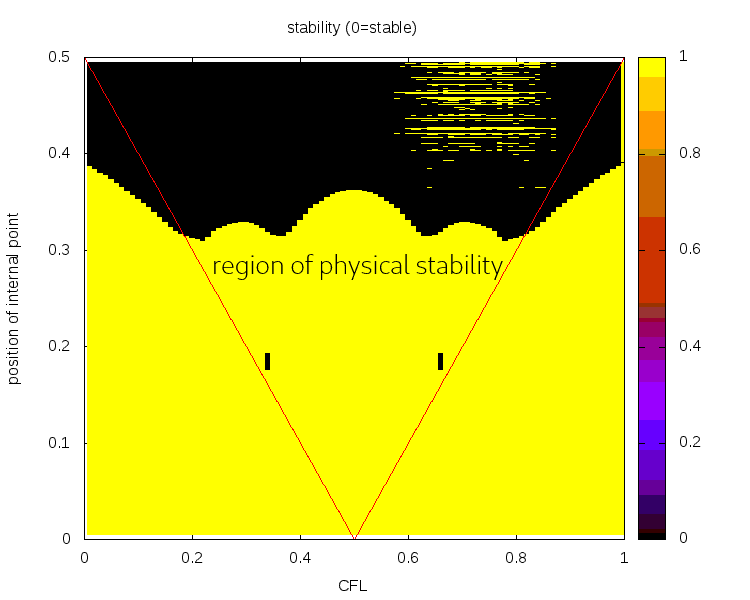}  \hfill \includegraphics[width=0.32\textwidth]{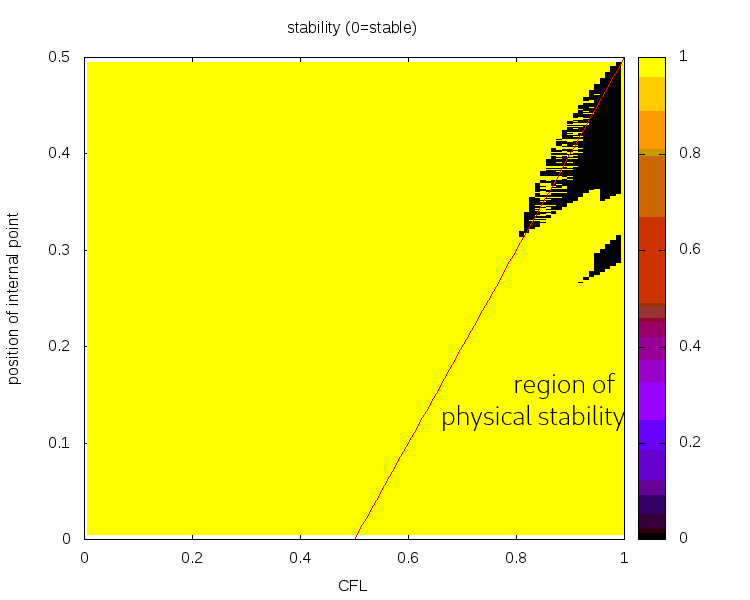}
 \caption{Two point values located symmetrically $\xi_2 = -\xi_1$ (fifth order method) with $M=5$ equidistant quadrature points in time. Lower panel shows stability as a function of the location of the right internal point value ($\xi_2$) and the CFL number (black = stable). The noise in the central and right figure has to do with the sampling and with one of the eigenvalues being exactly one; physical stability regions (i.e. those combinations of point location and CFL for which the foot of the characteristic is actually in the cell whose reconstruction is used in the update) are to the left, between and to the right of the red lines, respectively. One observes that actual stability also occurs outside the region of physical stability, while certain parameter choices inside the physical stability region turn out to be actually unstable.}
 \label{fig:points5thstability}
\end{figure}

Note also that according to the definition of $q_\text{recon}$ in \eqref{eq:pointsglobalrecon} the foot of the characteristic in \eqref{eq:pointsupdateintermedaite} is evaluated in the ``correct'' cell: if 
\begin{align}
 \Delta x \xi_j - c \Delta t > -\frac{\Delta x}{2}
\end{align}
then information is taken from cell $i$, otherwise, for
\begin{align}
 \Delta x \xi_j - c \Delta t \leq -\frac{\Delta x}{2}
\end{align}
from cell $i-1$, (see Figure \ref{fig:points5thstability}). 

\begin{figure}  % 2b3d6c2
 \centering
 \includegraphics[width=0.45\textwidth]{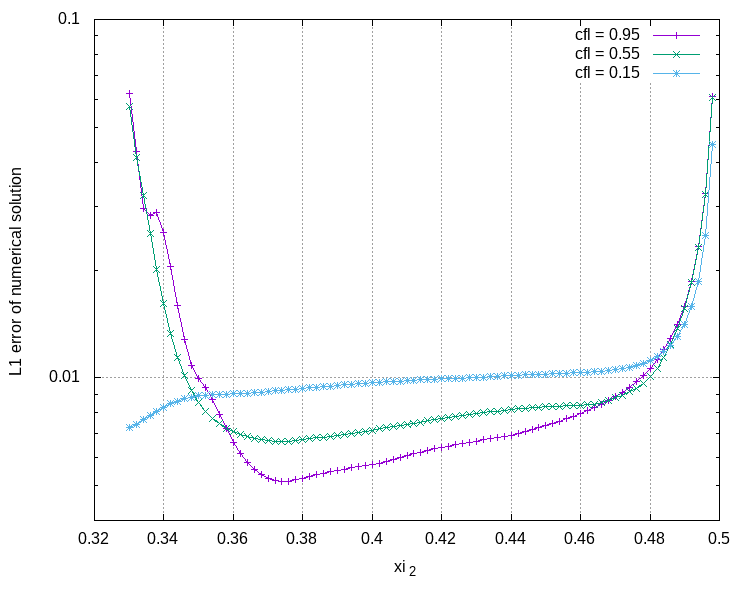} \hfill \includegraphics[width=0.45\textwidth]{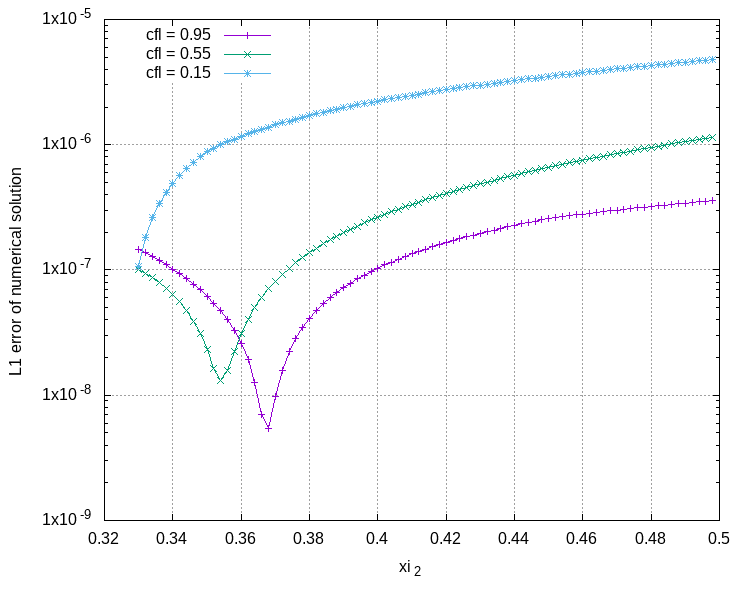}
 \caption{Dependence of the error on the location of the internal point values for different values of the CFL number for a 5th order method. Linear advection with speed $c=1$ is solved on a domain $[0,1]$ with periodic boundaries, discretized by 200 cells. The error of the point values is shown at time $t = 30$; the errors of the cell averages are very similar and are not shown. The three curves (for CFL numbers 0.15, 0.55 and 0.95) show the dependence of the error on the location $\xi_2$ of the right internal point value, with $\xi_1 = -\xi_2$. This corresponds to a 5th order method, and the fluxes are obtained through quadrature in time with $M=4$ Gauss-Lobatto points. As can be seen from Figure \ref{fig:points5thstability}, stability can be guaranteed for all CFL numbers only if $\xi_2 \gtrsim 0.365$, but for a CFL number of 0.15, $\xi_2 \simeq 0.33$ is still stable. This explains the low error of the corresponding curve for small $\xi_2$. \textit{Left}: The initial data are discontinuous with $q_0 = 1$ for $x \in [\frac14, \frac34]$, and $1/10$ elsewhere. \textit{Right}: The initial data are smooth with $q_0(x) = 0.8 + \exp(-(x-0.5)^2/0.05^2)$.}
 \label{fig:xierror}
\end{figure}

This guarantees optimal stability, as can be seen in Figure \ref{fig:points5thstability} for 5th order, which shows a stability analysis associated to different cases: for small CFL numbers, both characteristics used in the update of the internal point values have their foot points in the same cell; for intermediate CFL numbers, one of them reaches out into the neighbouring cell, and for large CFL numbers, the foot points of both are in the neighbouring cell. (I.e. the CFL number is based on $\Delta x$ and the largest characteristic speed, not on the smallest distance between point values.)
One observes that taken together, these yield a stability result valid for all CFL numbers between 0 and 1. At the same time, one remarks that, surprisingly, the location of the internal point values is not arbitrary. For stability, the internal point values have to be distributed fairly close to the cell interfaces, with $\xi_2 = - \xi_1 \gtrsim 0.365$. We are not aware of a simple reason that would explain this finding. 
At the same time, from Fig. \ref{fig:points5thstability} (middle) one deduces that the case of only one of the characteristics reaching out into the left cell is stable even outside the physical region of stability. Once inside the stable domain, we do not find the precise location of the point values to have a visible influence on the simulation results, see Figure \ref{fig:xierror}. One observes that the error grows as $\xi_2 \to \frac12$, and that it drops towards the stability limit. There does not seem to exist a unique value of $\xi_2$ that would dramatically enhance the properties of the method, but smaller values seem more accurate, as long as they respect the stability limit. The same analysis yields very similar results for Burgers' equation (not shown).

\subsection{Transonic upwinding} \label{ssec:transonicupwindingpoints}

We use the same algorithm as the one described in \cite{barsukow19activeflux} and in Section \ref{ssec:transonicupwindinggeneral} to which the reader is referred for more details. In short, for scalar conservation laws, the iterative approximation to the foot point of the characteristic passing through $x_{i+\frac12}$ at time $t$ is initialized with two values ($q^n_{i-\frac12}$ and $q^n_{i+\frac32}$ instead of $q^n_{i+\frac12}$), see Equation \eqref{eq:fixpointiterationinit}. Every iteration increases the order of accuracy of the evolution operator by one. After a sufficient number of iterations, the speeds of the two estimated characteristics are compared, and the one with the largest speed in absolute value is selected. This treatment of the transonic case does not affect the accuracy of the evolution operator. For systems, the approximate evolution operator is evaluated at different initialization locations. Then, the one with the largest sum of absolute values of the eigenvalues is chosen; for further details the reader is again referred to \cite{barsukow19activeflux}.

\subsection{Limiting}

The limiting strategy is to check whether the reconstruction exceeds the minimum/maximum of $\{ \bar q_i, q_{i\pm\frac12} \} \cup \{q_{i,j}\}_{j=1,\ldots,k}$, and to use a parabolic/power law reconstruction in that case. As a consequence, the internal point values are ignored and only the point values at cell interfaces contribute. In principle, it is possible to do this successively and to decrease the order by two until third order is reached. The decision whether to use a parabolic or a power-law reconstruction remains as described in Section \ref{ssec:limiting}. As it is too costly to evaluate the presence of local extrema analytically for high-degree polynomials, the value of the reconstruction is checked at 10 locations inside the cell.

\subsection{Numerical results}

\subsubsection{Linear advection}

\begin{figure} % 
  \centering
  \includegraphics[width=0.49\textwidth]{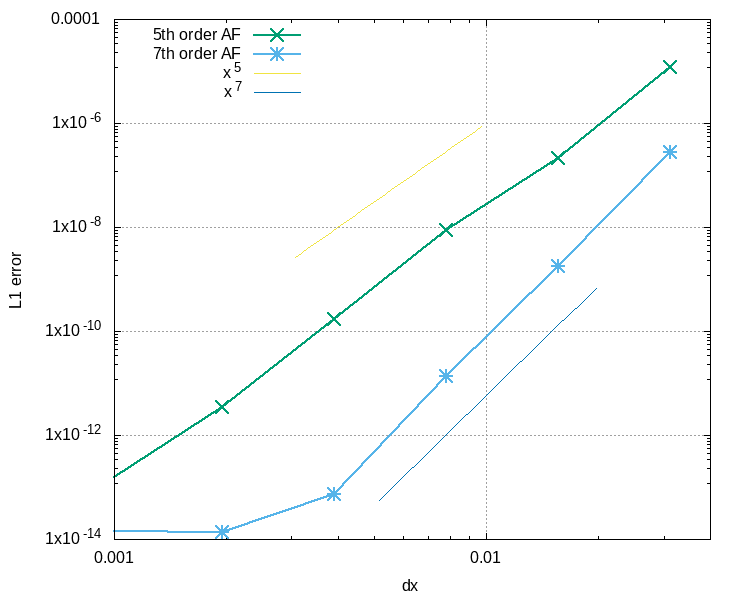}\hfill\includegraphics[width=0.49\textwidth]{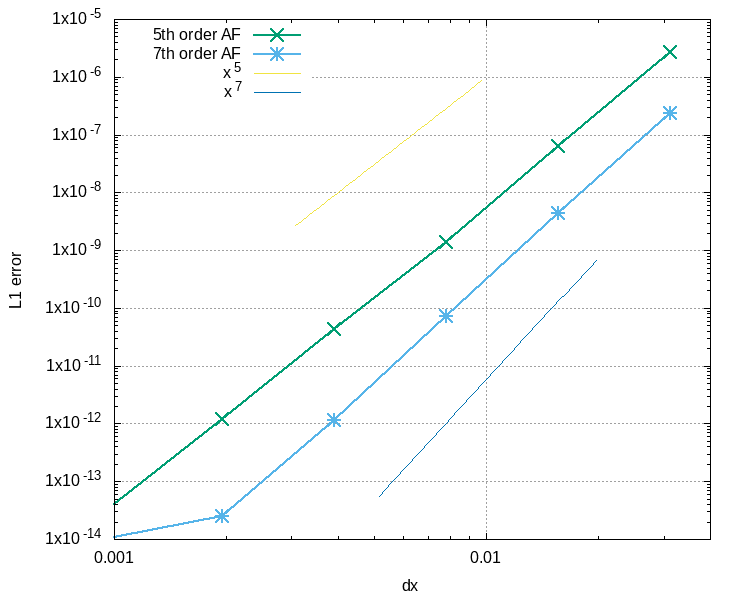}
  \caption{Convergence of the the high order extension of Active Flux via additional point values for linear advection. \textit{Left}: $L^1$ error of the point values. \textit{Right}: $L^1$ error of the averages.}
  \label{fig:points-advection-convergence}
 \end{figure}

For linear advection with speed $1$, the initial data \eqref{eq:initialconvergence} are evolved using a CFL number of 0.5 on grids covering $[0,1]$. For 5th order accuracy, two internal point values located at $\xi_{1,2} = \pm0.415$ are used, as well as a Gauss-Lobatto quadrature for \eqref{eq:pointsfluxquadrature} with $M=4$ points in time. For 7th order, 4 internal points are located at $\{\pm0.48, \pm0.41 \}$, and a Gauss-Lobatto quadrature with $M=5$ points in time is used. The exact evolution operator \eqref{eq:updatelinearadvectionvanleer} is used. The setup is run without limiting. Figure \ref{fig:points-advection-convergence} shows the experimental convergence at time $t=0.01$.

\subsubsection{Burgers' equation}

\begin{figure} % 7f0a164, ca30a6f
  \centering
\includegraphics[width=0.49\textwidth]{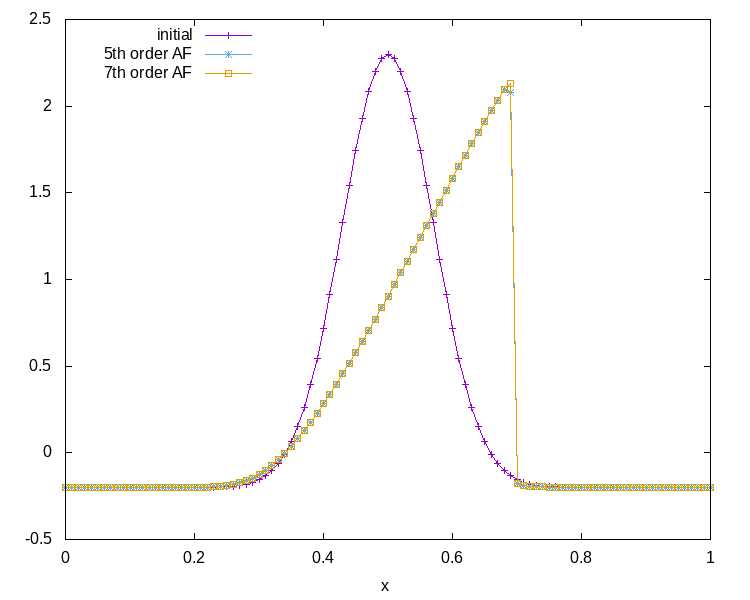}\hfill\includegraphics[width=0.49\textwidth]{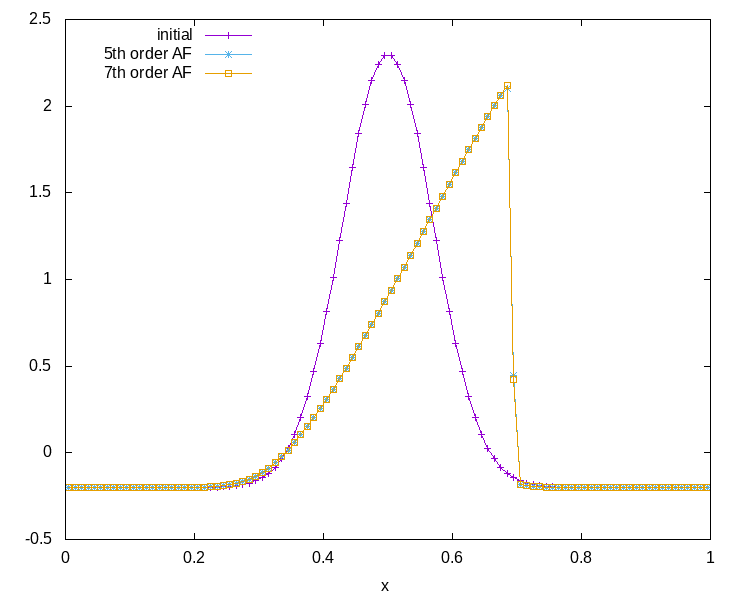}
  \caption{Numerical evolution of Gaussian initial data for Burgers' equation at time $t=0.1$. \textit{Left}: Point values. \textit{Right}: Averages.}
  \label{fig:points-burgers}
 \end{figure}

For Burgers' equation we consider again transonic initial data \eqref{eq:initialtransonicburgers} on a grid of 100 cells covering $[0,1]$. Figure \ref{fig:points-burgers} shows the numerical solution at time $t=0.1$ after shock formation. A CFL number of 0.5 has been used. For 5th order accuracy, two internal points values are located at $\xi_{1,2} = \pm0.415$ and a Gauss-Lobatto quadrature for \eqref{eq:pointsfluxquadrature} is used with $M=4$ points in time. For 7th order, 4 internal points are located at $\{ \pm0.48, \pm0.41\}$, and a Gauss-Lobatto quadrature is used with $M=5$ points in time. In both cases, the approximate evolution operator \eqref{eq:fixpointiterationinit}--\eqref{eq:fixpointiteration} is used, with an adequate number of fixpoint iterations according to \eqref{eq:fixpointiterationerror}. 

\begin{figure} %e0d6f19, 0a704cc
 \centering
  \includegraphics[width=0.49\textwidth]{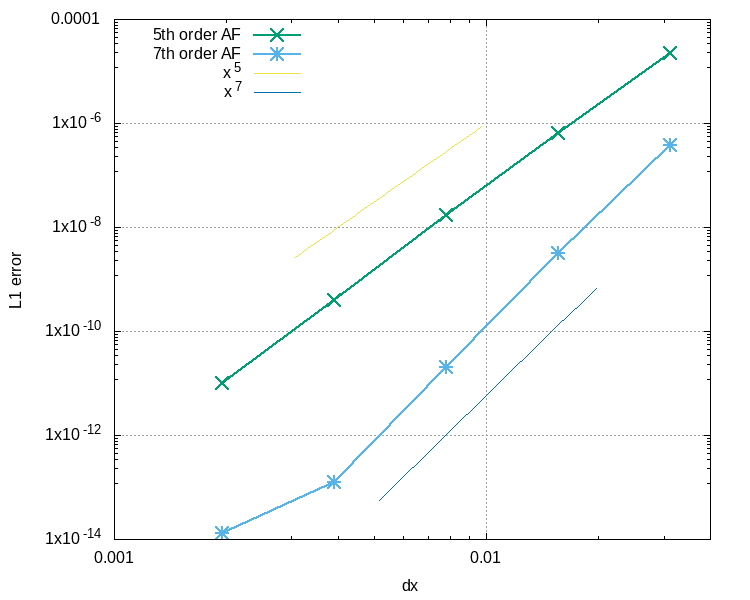}\hfill\includegraphics[width=0.49\textwidth]{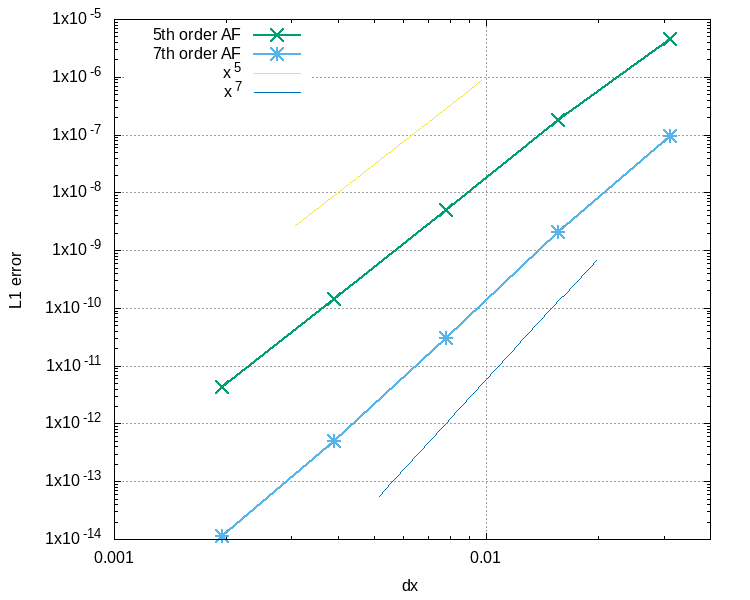}
  \caption{Convergence of the the high order extension of Active Flux via additional point values for Burgers' equation. \textit{Left}: $L^1$ error of the point values. \textit{Right}: $L^1$ error of the averages.}
  \label{fig:points-burgers-convergence}
\end{figure}

Figure \ref{fig:points-burgers-convergence} shows convergence results at time $t=0.01$ prior to shock formation. Here, limiting is not used, but the transonic upwinding follows the procedure described in \ref{ssec:transonicupwindingpoints}.

\subsubsection{Euler equations}

\begin{figure} % f8afdfb
 \centering
 \includegraphics[width=0.49\textwidth]{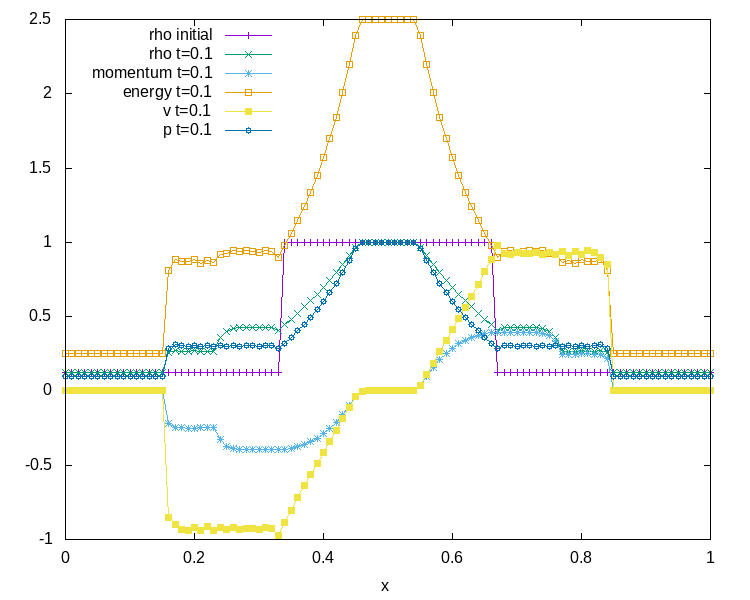}\hfill\includegraphics[width=0.49\textwidth]{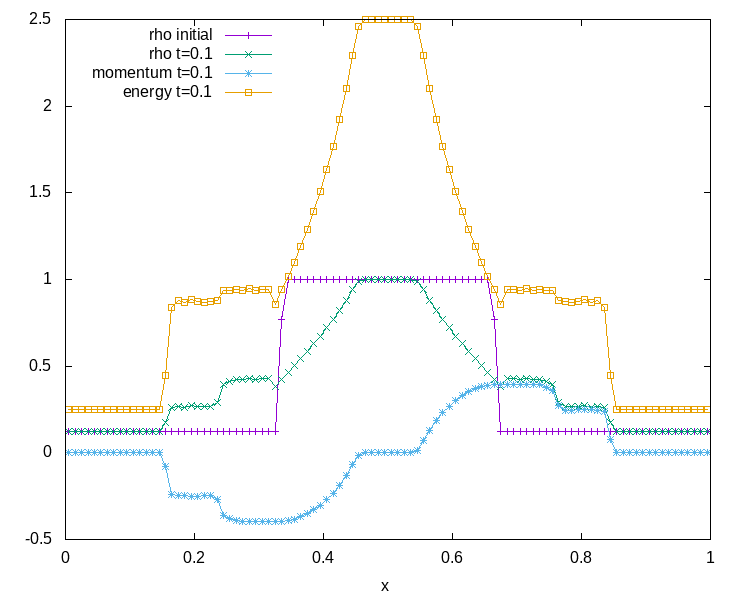}
  \caption{Numerical evolution of the Sod shock tube at time $t=0.1$ using the high order extension of Active Flux via additional point values. \textit{Left}: Point values. \textit{Right}: Averages.} 
 \label{fig:pointssod}
\end{figure}

The Sod shock tube is evolved on a grid of 100 points covering $[0,1]$, using a CFL number of 0.5. The additional point values are located at $\xi_{1,2} = \pm 0.415$, the flux quadrature \eqref{eq:pointsfluxquadrature} is performed using a Gauss-Lobatto rule involving $M=4$ points in time, a setup corresponding to 5th order accuracy. The approximate evolution operator is described in \cite{barsukow19activeflux}. It is formally third order accurate only; however, this can be seen in parallel to the common usage of low order time integrators (e.g. RK3) for spatially high order methods. An extension of the approximate evolution operator from \cite{barsukow19activeflux} to higher order does not seem generally impossible, but for now remains subject of future work. Numerical results at $t=0.1$ are shown in Figure \ref{fig:pointssod}.

\section{High order via further moments} \label{sec:moments}

Finally, it is possible to interpret Active Flux in a way closer to Finite Element methods, by declaring higher moments to be new degrees of freedom in a high-order version. It is thus a complementary view to the previous section. 

Define the $p$-th moment as
\begin{align}
 q_i^{(p)} := A_p \int_{-\frac{\Delta x}{2}}^{\frac{\Delta x}{2}} \dd x \, x^p q(t,x_i + x) 
\end{align}
The normalization $A_p = \frac{(p+1) 2^p}{\Delta x^{p+1}}$ is chosen such that $q(t,x) \equiv 1$ implies $q^{(p)} = 1$ $\forall p \in 2\mathbb N_0$. The first values are $(A_p)_p = \left( \frac{1}{\Delta x}, \frac{4}{\Delta x^2}, \frac{12}{ \Delta x^3} ,\frac{32}{\Delta x^4} , \frac{80}{\Delta x^5} , \ldots \right)$ (sequence A001787 in the Online Encyclopedia of Integer Sequences (OEIS)).

With $k$ moments, every cell has access to $k + 2$ pieces of information and thus, a polynomial reconstruction of degree $k+1$ can be constructed in every cell. The order of accuracy of the method is thus expected to be $k+2$. As the point values are shared, on average there are $k+1$ degrees of freedom per cell. The global reconstruction remains $C^0$.

\subsection{Update of the moments}

Multiplying \eqref{eq:conslaw} by $x^p$ and integrating yields the evolution equation for the moments:
\begin{align}
 A_p \int_{-\frac{\Delta x}{2}}^{\frac{\Delta x}{2}} \dd x \, x^p \del_t q  + A_p \int_{-\frac{\Delta x}{2}}^{\frac{\Delta x}{2}} \dd x \, x^p \del_x f(q)  &= 0\\
 \frac{\dd}{\dd t} q^{(p)}_i + A_p \left( \left( \frac{\Delta x}{2}\right)^p f(q_{i+\frac12}) - \left( -\frac{\Delta x}{2}\right)^p f(q_{i-\frac12}) \right) -p A_p \int_{-\frac{\Delta x}{2}}^{\frac{\Delta x}{2}} \dd x \, x^{p-1} f(q)  &= 0\\
 \frac{\dd}{\dd t} q^{(p)}_i + (p+1)  \frac{ f(q_{i+\frac12}) - ( -1)^p f(q_{i-\frac12})  }{\Delta x}
 -p A_p \int_{-\frac{\Delta x}{2}}^{\frac{\Delta x}{2}} \dd x \, x^{p-1} f(q) &= 0 \label{eq:advectionmomentsgeneral}
\end{align}
For $p=0$ one just has the usual equation for the average $\bar q_i \equiv q_i^{(0)}$
\begin{align}
 \frac{\dd}{\dd t} q^{(0)}_i +  \frac{ f(q_{i+\frac12}) - f(q_{i-\frac12})  }{\Delta x} &= 0 \label{eq:advectionmoments0}
\end{align}

For linear flux $f(q) = cq$, the integral in \eqref{eq:advectionmomentsgeneral} can be expressed as a moment of one order less:
\begin{align}
\frac{\dd}{\dd t} q^{(p)}_i + (p+1) c \frac{ q_{i+\frac12} - ( -1)^p q_{i-\frac12} }{\Delta x}
 -\frac{2(p+1)}{\Delta x} \frac{p 2^{p-1}}{\Delta x^p} c \int_{-\frac{\Delta x}{2}}^{\frac{\Delta x}{2}} \dd x \, x^{p-1} q &= 0\\
\frac{\dd}{\dd t} q^{(p)}_i + (p+1) c \frac{ q_{i+\frac12} - ( -1)^p q_{i-\frac12} }{\Delta x}
 - \frac{2(p+1)}{\Delta x}  c q^{(p-1)}_i &= 0 \quad (p\geq 1) \label{eq:advectionmoments}
\end{align}
such that, for example, the first three moments thus evolve according to
\begin{align}
\frac{\dd}{\dd t} q^{(0)}_i + c \frac{ q_{i+\frac12} - q_{i-\frac12} }{\Delta x}  &= 0\\
\frac{\dd}{\dd t} q^{(1)}_i + 2 c \frac{ q_{i+\frac12} + q_{i-\frac12} }{\Delta x} - \frac{4}{\Delta x}  c q^{(0)}_i &= 0\\
\frac{\dd}{\dd t} q^{(2)}_i + 3 c \frac{ q_{i+\frac12} - q_{i-\frac12} }{\Delta x} - \frac{6}{\Delta x}  c q^{(1)}_i &= 0
\end{align}

Observe that the evolution of the moments is exact (up to a possible quadrature error upon evaluating the integral in \eqref{eq:advectionmomentsgeneral} in case of a nonlinear $f$). Apart from quadrature errors, the accuracy of the method is determined entirely by the choice of the point value update. 

\begin{definition}[Active Flux with higher moments]
The \emph{Active Flux method with higher moments} is the following semi-discretization
\begin{align}
 \begin{cases} \displaystyle  \frac{\dd}{\dd t} q^{(p)}_i = - (p+1)  \frac{ f(q_{i+\frac12}) - ( -1)^p f(q_{i-\frac12})  }{\Delta x}
 +p A_p \int_{-\frac{\Delta x}{2}}^{\frac{\Delta x}{2}} \dd x \, x^{p-1} f(q)  \quad p = 0, \ldots, k \\ \\
 \displaystyle \frac{\dd}{\dd t} q_{i+\frac12}(t) = -R\Big(q_{i-\frac12}(t),  q_{i}^{(0)}(t), \ldots, q_i^{(k)}(t), q_{i+\frac12}(t), q_{i+1}^{(0)}(t), \ldots, q_{i+1}^{(k)}(t), q_{i+\frac32}(t)\Big) \\\nonumber \phantom{mmmmmmmmmm} i \in \mathbb Z, k\geq 0
 %- \Big (  f'(\tilde q_{i+\frac12})^+ D  + f'(\tilde q_{i+\frac12})^- D^* \Big ) \label{eq:pointvalueupdatefindiff}\qquad i \in \mathbb Z, k\geq 0, m \geq 0 
 \end{cases}
\end{align}
of \eqref{eq:conslaw} with $A_p = \frac{(p+1)2^p}{\Delta x^{p+1}}$ and the interpretations
\begin{align}
 q_i^{(p)}(t) &\simeq A_p \int_{x_{i-\frac12}}^{x_{i+\frac12}} \dd x \, x^p q(t, x)  &
 q_{i+\frac12}(t) &\simeq q(t, x_{i+\frac12})
\end{align}
and $R$ a consistent approximation of $\del_x f(q)$ at $x_{i+\frac12}$.
\end{definition}

\subsection{Update of the point values}

We propose to update the point value with a finite difference approximation to the derivative, obtained as the derivative of the reconstruction. This means that this finite difference has a compact stencil, and only involves the degrees of freedom of one of the adjacent cells. 

For an $N$-th order accurate method, define the unique polynomial $q_{\text{recon},i}: \left[ -\frac{\Delta x}{2}, \frac{\Delta x}{2} \right] \to \mathbb R^m$ of order $N-1$ fulfilling
\begin{align}
 q_{\text{recon},i}\left( \pm \frac{\Delta x}{2} \right ) &= q_{i\pm\frac12} &
  A_p \int_{-\frac{\Delta x}{2}}^{\frac{\Delta x}{2}} \dd x \, x^p q_{{recon},i}(x)  &= q_i^{(p)} \qquad p = 0, \ldots, N-3
\end{align}

%\todo{does this polynomial always exist?}

For example, the lowest order polynomial with given values $q_{i\pm\frac12}$ and moments 0 through 2 is
\begin{align}
 \frac{1}{16} q_{i+\frac12} (2\xi+1) (3 - 30 \xi - 60 \xi^2 + 280 \xi^3) + \frac{1}{16} q_{i-\frac12} (2\xi-1) (-3 - 30 \xi + 60 \xi^2 + 280 \xi^3) \label{eq:momentspolynomialexample}\\\nonumber+ (2\xi-1) (2\xi+1) \left( \frac{15}{16} q^{(0)}  (28 \xi^2 - 3) - \frac{15}2 q^{(1)} \xi  - \frac{35}{16} q^{(2)}  (20 \xi^2 - 1) \right ) 
\end{align}
where $\xi = x/\Delta x$.

A finite difference approximation to the derivative at $x=\frac{\Delta x}{2}$ is obtained by evaluating $(\del_x q_{\text{recon},i})\left( \frac{\Delta x}{2} \right )$. It necessarily is a linear function of $q_{i+\frac12}, q_{i-\frac12}, q^{(0)}_i , \ldots, q^{(N-3)}_i$ and therefore can be written as
\begin{align}
 \frac{1}{\Delta x} \left(  b_{\frac12} q_{i+\frac12} + b_{-\frac12} q_{i-\frac12} + \sum_{j=0}^{N-3} b_0^{(j)} q^{(j)}_i   \right) = q'(x_{i+\frac12}) + \mathcal O(\Delta x^N) \label{eq:momentdifference}
\end{align}

\begin{notation}
 A finite difference approximation of the form \eqref{eq:momentdifference} shall be depicted graphically as follows:
\begin{align}
 \frac{1}{\Delta x} \left(  b_{\frac12} q_{i+\frac12} + b_{-\frac12} q_{i-\frac12} + \sum_{j=0}^{N-3} b_0^{(j)} q^{(j)}_i   \right) = 
 \begin{array}{|c|ccc||c|}
   \hline &&&&\\[-12pt]
   && b_0^{(N-3)} && \\
   && \vdots && \\
   && b_0^{(1)} && \\
   && b_0^{(0)} && \\
   \multicolumn{2}{|c}{b_{-\frac12}} && \multicolumn{2}{c|}{b_{\frac12}} \\[4pt]\hline
 \end{array}
\end{align}
Observe the consistency of the notation with that introduced in \eqref{eq:definitiontableau}. Note that, again, the double vertical line indicates the cell interface at which the finite difference formula is an approximation of the derivative to the given order.

\end{notation}

The finite differences constructed in the way described above are unique for every order of accuracy; to lowest orders they read\footnote{MD is supposed to mean \emph{moment difference}.}

\begin{align}
 \text{MD}3 &= 
 \begin{array}{|c|ccc||c|}
   \hline &&&&\\[-12pt]
   && -6 && \\
   \multicolumn{2}{|c}{2} && \multicolumn{2}{c|}{4} \\[4pt]\hline
 \end{array} &
 \text{MD}5 &= \begin{array}{|c|ccc||c|}
   \hline &&&&\\[-12pt]
   && -35 && \\
   && -15 && \\
   && 15 && \\
   \multicolumn{2}{|c}{4} && \multicolumn{2}{c|}{16} \\[4pt]\hline
 \end{array} &
 \text{MD}7 &= \begin{array}{|c|ccc||c|}
   \hline &&&&\\[-12pt]
   && -\frac{693}4 && \\[2pt]
   && -\frac{315}4  && \\[2pt]
   && \frac{315}2 && \\[2pt]
   && \frac{105}2 && \\[2pt]
   && -\frac{105}4 && \\[2pt]
   \multicolumn{2}{|c}{6} && \multicolumn{2}{c|}{36} \\[4pt]\hline
 \end{array} \label{eq:momentdifferenceexamples}
\end{align}

Note that the coefficients of even moments and point values have to add up to 0. For linear advection $\del_t q + c \del_x q =0$ for positive $c$, the update with MD5 amounts to
\begin{align}
\del_t q_{i+\frac12} = - c \frac{4 q_{i-\frac12} + 15 q_i^{(0)} - 15 q_i^{(1)} - 35 q_i^{(2)} + 16 q_{i+\frac12} }{\Delta x}
\end{align}

An Active Flux method of order $p$ is obtained by using a finite difference MD$p$, just as is the case for finite differences in Section \ref{sec:findiff}. For a seventh-order method, the derivative approximation is MD7
\begin{align}
 3\frac{8 q_{i-\frac12} -35 q_i^{(0)} + 70  q_i^{(1)} +210 q_i^{(2)} - 105 q_i^{(3)} - 231 q_i^{(4)} + 48 q_{i+\frac12}}{4\Delta x}
\end{align}

Define again by $D^*$ the flipped version of the finite difference $D$ as follows:

\begin{align}
 D &= \begin{array}{|c|ccc||c|}
   \hline &&&&\\[-12pt]
   && b_0^{(N-3)} && \\
   && \vdots && \\
   && b_0^{(2)} && \\
   && b_0^{(1)} && \\
   && b_0^{(0)} && \\
   \multicolumn{2}{|c}{b_{-\frac12}} && \multicolumn{2}{c|}{b_{\frac12}} \\[4pt]\hline
 \end{array} &
 D^* &= \begin{array}{|c||ccc|c|}
   \hline &&&&\\[-12pt]
   && \pm b_0^{(N-3)} && \\
   && \vdots && \\
   && -b_0^{(2)} && \\
   && b_0^{(1)} && \\
   && -b_0^{(0)} &\phantom{n}& \\
   \multicolumn{2}{|c}{-b_{\frac12}} && \multicolumn{2}{c|}{-b_{-\frac12}} \\[4pt]\hline
 \end{array}
\end{align}

Observe that upon a reflection of space $x \mapsto -x$, the sign of the derivative changes, but also that of odd moments.

These flipped finite differences can of course be obtained by evaluating the derivative of the polynomial reconstruction at $-\frac{\Delta x}{2}$ and by performing an index shift. For example, the derivative of \eqref{eq:momentspolynomialexample} at $\frac{\Delta x}{2}$ is
\begin{align}
 \frac{ 4 q_{i-\frac12} + 15 q_i^{(0)} - 15 q_i^{(1)} - 35 q_i^{(2)} + 16 q_{i+\frac12}}{\Delta x}
\end{align}
while the one at $-\frac{\Delta x}{2}$ is
\begin{align}
 \frac{-4 q_{i+\frac12}  - 15 q_i^{(0)} - 15 q_i^{(1)} + 35 q_i^{(2)} - 16 q_{i-\frac12}}{\Delta x}
\end{align}

Given these finite differences, the semi-discrete method is Equations \eqref{eq:advectionmomentsgeneral} and
%\begin{align}
% \del_t q_{i+\frac12}(t) = - c (\del_x q_{\text{recon},i})\left( \frac{\Delta x}{2} \right )
%\end{align}
\begin{align}
 \frac{\dd}{\dd t} q_{i+\frac12}(t) = - \Big (  f'(\tilde q_{i+\frac12})^+ D  + f'(\tilde q_{i+\frac12})^- D^* \Big ) \label{eq:updatepointvaluemoments}
\end{align}
The point value update is thus in close analogy to \eqref{eq:pointvalueupdatefindiff}, with $D$ now being one of the moment differences \eqref{eq:momentdifferenceexamples}, and the method using MD3 is even exactly the same as \eqref{eq:pointvalueupdatefindiff} with FD3.

\subsection{Time integration and stability results}

Equations \eqref{eq:advectionmomentsgeneral} and \eqref{eq:updatepointvaluemoments} constitute a system of ODEs that we propose to integrate in time using e.g. SSP-RK3. For linear advection, MD5 is found to be stable for $\text{CFL} \leq 0.13 $ when using RK3 (see Section \ref{app:stability}). Experiments also suggest stability $\text{CFL} < 0.17$ for RK5. For MD7, RK3 yields stability until $\text{CFL} \leq 0.066$ (experimentally, $\text{CFL} \leq 0.1$ for RK7).

\subsection{Limiting}

As in the previous versions, we propose to gradually decrease the order of the method upon detection of oscillations. The finite differences employed in the point value update are constructed by differentiating\footnote{In fact, for the ease of the algorithm we perform the differentiation numerically.} a reconstruction. If the reconstruction is found to be non-monotone, the highest moment is neglected and a reconstruction of lower polynomial degree considered, until one either arrives at a monotone reconstruction, or at the parabolic reconstruction that employs the two point values and the average. In case it also is non-monotone, we propose once more to follow the limiting strategy of \cite{barsukow19activeflux} by replacing the parabola by a power law \eqref{eq:limitingpowerlaw1}--\eqref{eq:limitingpowerlaw2}, whenever a monotone reconstruction is at all possible. Their derivatives at cell interfaces are given in \eqref{eq:powerlawderivative1}--\eqref{eq:powerlawderivative2}. For polynomials of degree higher than 2, monotonicity is checked approximately upon evaluation of the polynomial at 5 locations inside the cell.

\subsection{Transonic upwinding}

We employ the same choice of $\tilde q_{i+\frac12}$ as in Section \ref{ssec:findifftransonicupwinding}.

\subsection{Numerical results}

\subsubsection{Linear advection}

\begin{figure} % f3c9b3c
  \centering
  \includegraphics[width=0.49\textwidth]{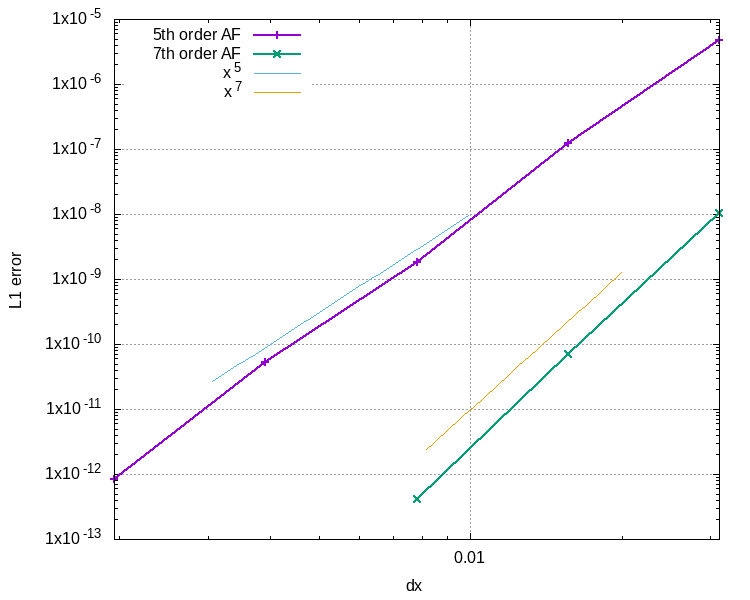}\hfill\includegraphics[width=0.49\textwidth]{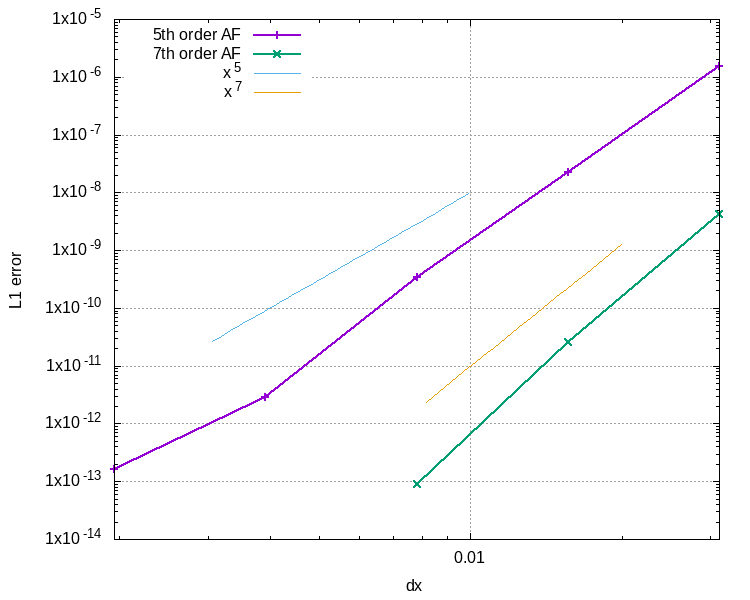}
  \caption{Convergence of the the high order extension of Active Flux via higher moments for linear advection. \textit{Left}: $L^1$ error of the point values. \textit{Right}: $L^1$ error of the averages.}
  \label{fig:moments-advection-convergence}
 \end{figure}

Convergence results for linear advection with speed 1 are shown in Figure \ref{fig:moments-advection-convergence} at time $t=0.01$. Here, limiting is not used, and as before a very small CFL number of $10^{-4}$ is used with RK3 in order for the error to be dominated by that of spatial discretization. The initial data are again \eqref{eq:initialconvergence}.

\subsubsection{Burgers' equation}

\begin{figure} % debe7be
  \centering
\includegraphics[width=0.49\textwidth]{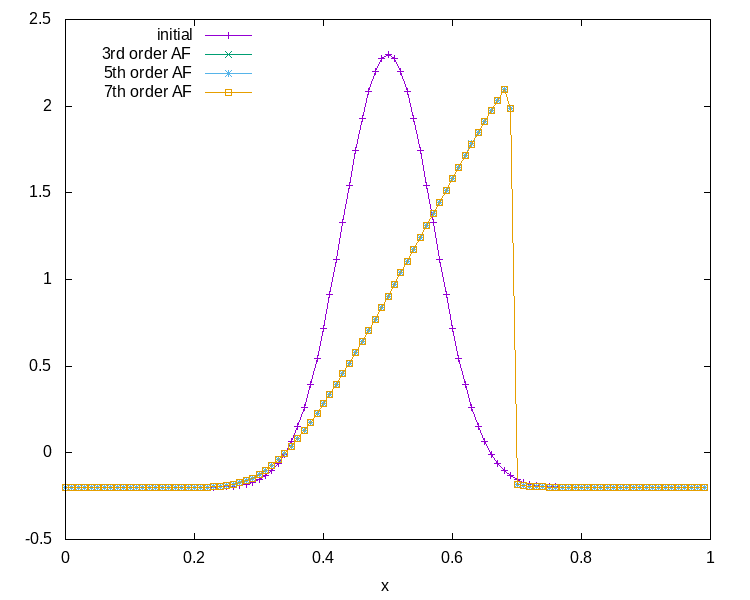}\hfill\includegraphics[width=0.49\textwidth]{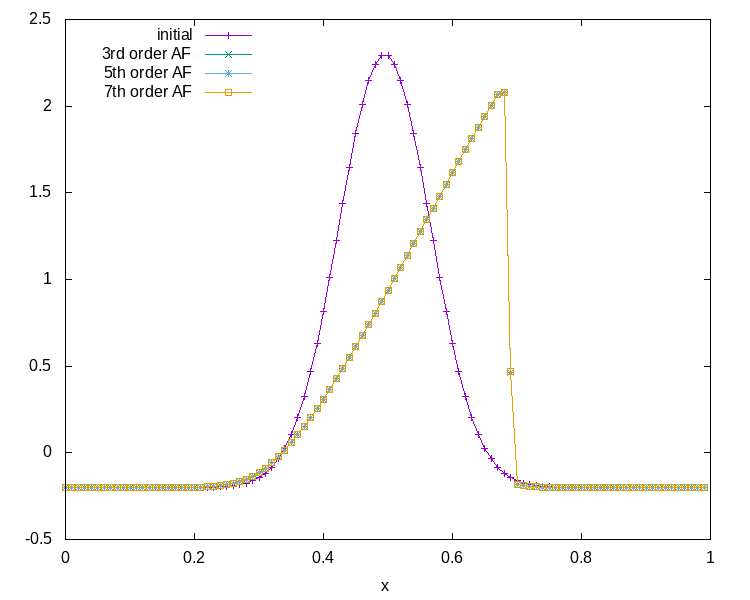}
  \caption{Numerical evolution of Gaussian initial data for Burgers' equation at time $t=0.1$. Setup as in Section \ref{ssec:findiffburgers}. \textit{Left}: Point values. \textit{Right}: Averages.}
  \label{fig:moments-burgers}
 \end{figure}

Figure \ref{fig:moments-burgers} shows the numerical evolution of transonic initial data \eqref{eq:initialtransonicburgers} for Burgers' equation for high order Active Flux via higher moments at $t=0.1$. The grid consists of 100 cells covering $[0,1]$ and the CFL number is 0.01. Limiting is used, as is the transonic upwinding modification. The shock is resolved well and no artefacts at the sonic point are visible.

\begin{figure} % 
  \centering
  \includegraphics[width=0.49\textwidth]{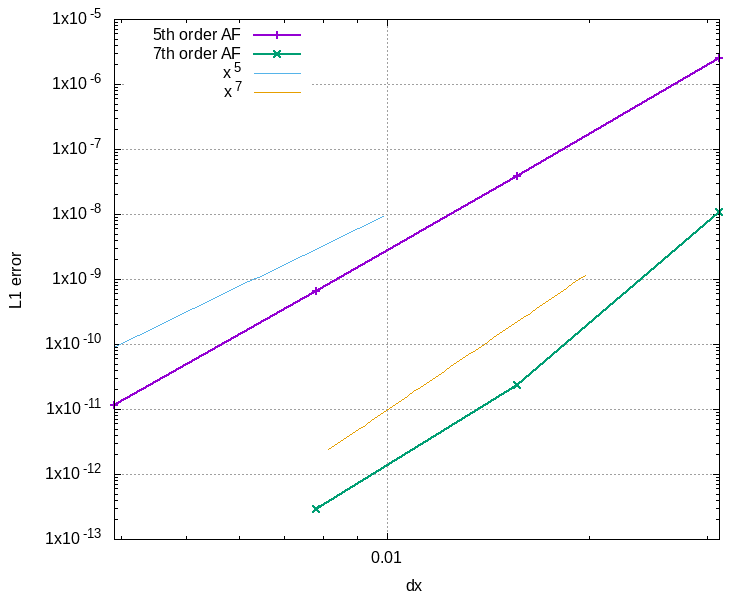}\hfill\includegraphics[width=0.49\textwidth]{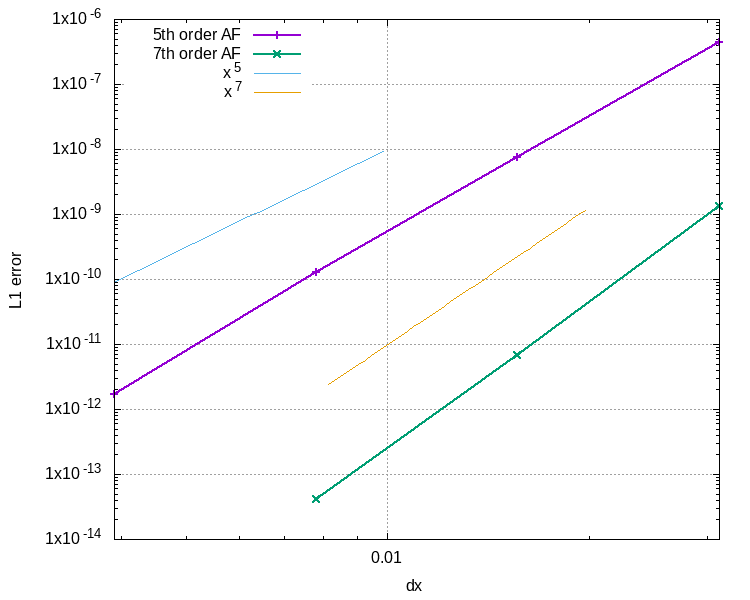}
  \caption{Convergence of the the high order extension of Active Flux via higher moments for Burgers' equation. \textit{Left}: $L^1$ error of the point values. \textit{Right}: $L^1$ error of the averages.}
  \label{fig:moments-burgers-convergence}
 \end{figure}

 Figure \ref{fig:moments-burgers-convergence} shows the experimental order of convergence for initial data \eqref{eq:initialconvergence} at time $t=0.01$, before the formation of a shock. No limiting is used, and the CFL number is again $10^{-4}$ (RK3) in order to ensure that the error is dominated by that of spatial discretization.
 
 %\todo{rerun}

\subsubsection{Euler equations}

Finally, a Sod shock tube is solved on a grid of 100 cells covering $[0,1]$. The initial discontinuities are located at $x=\frac13$ and $x=\frac23$, such that it is possible to use periodic boundaries. The numerical method employed is of 7th order accuracy in space; time integration is performed using RK3 with a CFL number of 0.01.

\begin{figure} % b2c9c9f
 \centering
 \includegraphics[width=0.49\textwidth]{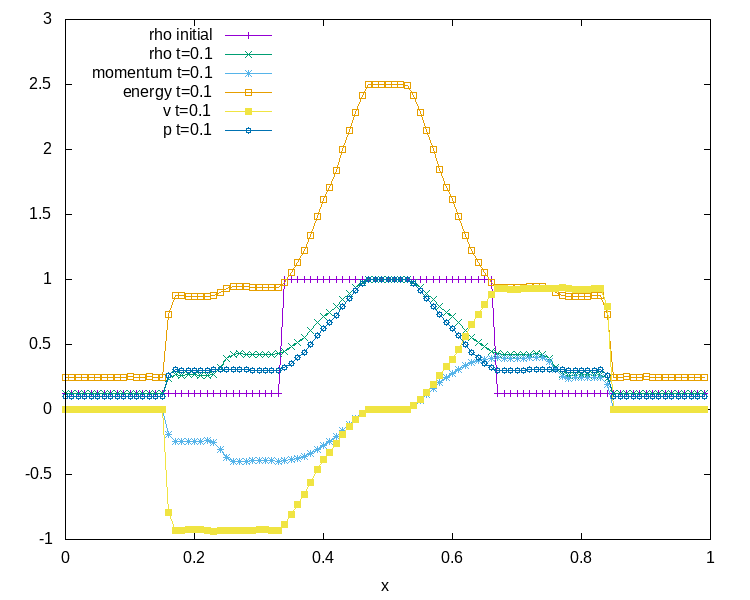}\hfill\includegraphics[width=0.49\textwidth]{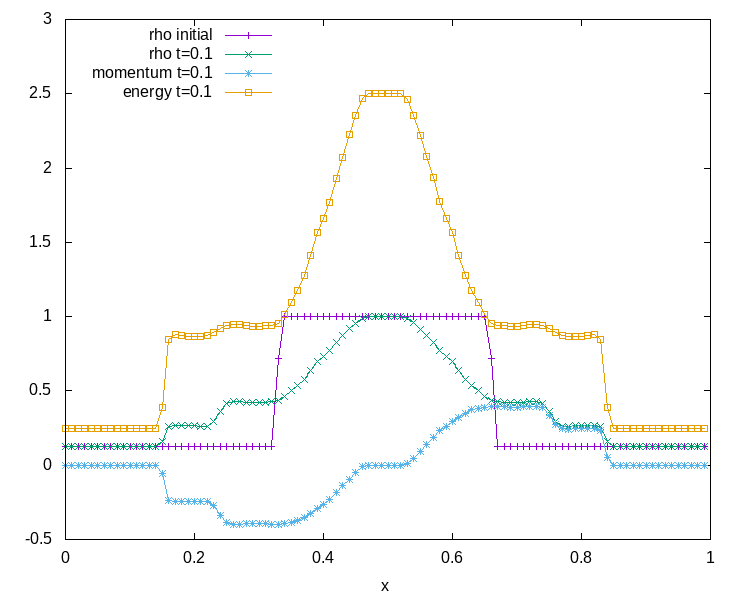}
  \caption{Numerical evolution of the Sod shock tube at time $t=0.1$. \textit{Left}: Point values. \textit{Right}: Averages.}
  \label{fig:momentseulersod7}
\end{figure}

\section{Higher order via higher derivatives}

Yet another high-order extension of Active Flux is presented in \cite{roe21}, where \emph{spatial derivatives} are stored at cell interfaces. When in every cell an average is stored, and at every cell interface a point value and derivatives of order $\leq s$, the modified Hermite interpolation problem gives rise to a polynomial $P^{2 s + 2}$ in every cell, and simultaneously augments the global regularity of the reconstruction to $C^s$. The average number of degrees of freedom per cell is only $s+2$ and the order of the scheme is no higher than $2 s + 3$. More information is stored at the cell interfaces, which reduces the memory  requirements of the method. Note, though, that the spatial derivatives cannot be used to increase the order of accuracy for the flux quadrature in time. Also, the evolution of spatial derivatives may be cumbersome, and the extension of the approach to multiple spatial dimensions at the moment remains unexplored. The approximate evaluation of the Schur-Miller criterion (Section \ref{app:stability}) suggests that, for linear advection, the 5th order version of this method is stable up to $\text{CFL} = 1$.

\section{Conclusions and outlook}

Active Flux, as it was presented so far since its first appearance in \cite{vanleer77}, was restricted to third order accuracy. During this time, the degrees of freedom of Active Flux in one spatial dimension always have been a cell average and a point value at every cell interface, which means that a parabola could naturally be used as a reconstruction in each cell. Third order accuracy has been so natural, that any attempt to go beyond third order must question the current foundations of the method, the understanding of ``what'' Active Flux is. In this paper, we associate three choices of an extension to higher order with three different interpretations of Active Flux. In all these interpretations we, however, have maintained two special properties of Active Flux, that we consider to be its distinctive features: continuity of the approximation and the presence of point values and cell averages among the degrees of freedom. The solution is considered globally continuous and, contrary to Godunov methods, Riemann Problems are never introduced in order to compute a numerical flux. The presence of point values (whose update includes upwinding) ensures stability, and allows to use non-conservative formulations of the equations for their update. At the same time, the presence of an average that is updated conservatively ensures the correct approximation of weak solutions.

The first view of Active Flux presented in this paper is that of a coupled Finite Volume/Finite Difference method (A) with a cell average updated according to the Finite Volume approach and point values, whose update includes finite difference approximations on a stencil whose size depends on the order of accuracy. The second view is that of an enriched Finite Volume method (B), with a cell average, and point values not only on cell interfaces but also inside the cell. The third view is that of a coupled Finite Element/Finite Difference method (C), with degrees of freedom being moments of the solution (starting with the average) and point values at cell interfaces that are used to update the moments.

The three high order extensions have very different properties. The characteristics-based update of the point values inside the cell for the variant B is the one that is most unorthodox, but also the one with the best stability properties. It is also closest to the time integration employed in the original Active Flux method. The usage of moments (variant C) as additional degrees of freedom maintains the elegance of a compact stencil and is the best candidate for further extensions, but also has the smallest CFL numbers among the considered variants. Both variants A and C come at the advantage that their application to any hyperbolic system is immediately available (the diagonalization of the Jacobian being the only necessary ingredient). 

Future work will be devoted to an extension of these methods to multiple spatial dimensions. As the Active Flux method is fairly recent (at least for nonlinear problems), certain questions that might seem settled for Finite Volume methods are still far from being conclusively answered for Active Flux. This concerns limiting as well as the convergence to the correct entropy solution. Although this paper presents viable approaches, more theoretical studies will be necessary, which surely will bring further improvements.

\bibliographystyle{alpha}

\appendix
\section{Stability analysis} \label{app:stability}

On equidistant Cartesian grids with $x_i = i \Delta x$, consider inserting a Fourier mode\footnote{$\ii$ is the imaginary unit.} $\hat q \exp(\ii K x)$ with $\hat q \in \mathbb C^m$, $K\in\mathbb R$ as data at time $t^n$. This means for point values
\begin{align}
 q_{i+\frac12}^n &= \hat q \exp\left( \ii K \frac{\Delta x}{2} \right ) \exp(\ii K i \Delta x) =: \hat q_p \exp(\ii K i \Delta x)\\
 q_{i+\frac32}^n &= \hat q_p \exp(\ii K \Delta x) \exp(\ii K i \Delta x)\\
 q_{i-\frac12}^n &= \hat q_p \exp(-\ii K \Delta x) \exp(\ii K i \Delta x)
\end{align}
Thus, every shift by one grid cell appears as an algebraic factor of $\exp(\ii K \Delta x) =: t_x$:
\begin{align}
 q_{i+\frac32} &= q_{i+\frac12} t_x & q_{i-\frac12} &= q_{i+\frac12} t_x^{-1} &&\text{etc.}
\end{align}

Similarly, for cell averages one finds
\begin{align}
 \bar q_i &= \hat q \frac{1}{\Delta x} \int_{x_i - \frac{\Delta x}{2}}^{x_i + \frac{\Delta x}{2}}\dd x \,  \exp(\ii K x) = \hat q \frac{1}{\Delta x} \exp(\ii K i \Delta x) \frac{\exp\left( \ii K \frac{\Delta x}{2} \right ) -\exp\left(- \ii K \frac{\Delta x}{2} \right )  }{\ii K}\\
 &= \hat q \,\sinc\left( K \frac{\Delta x}{2} \right )   \exp(\ii K i \Delta x) =: \hat q_a  \exp(\ii K i \Delta x)
\end{align}
with $\sinc x := \frac{\sin x}{x}$. Therefore, again,
\begin{align}
 \bar q_{i+1} &= \bar q_i t_x & \bar q_{i-1} &= \bar q_i t_x^{-1} &&\text{etc.}
\end{align}
The very same happens for moments.

Thus, any linear finite difference formula involving point values and moments can be written as a linear combination of just the point value and moments associated to one cell. For example, the update
\begin{align}
 \bar q_i^{n+1} &= \bar q_i^n - \frac{c\Delta t}{\Delta x}  (q_{i+\frac12} - q_{i-\frac12})
\end{align}
of $\bar q_i$ for linear advection upon inserting a Fourier mode reads
\begin{align}
 \bar q_i^{n+1} &= \bar q_i^n - \frac{c\Delta t}{\Delta x} q_{i+\frac12} (1 - t_x^{-1})
\end{align}

In general therefore, any of the three high-order versions of Active Flux presented in this paper can be written as
\begin{align}
 Q^{n+1} &= \mathcal A \, Q^n
\end{align}
with $A \in \mathbb C^{(2+N_1 + N_2) \times (2+N_1 + N_2)}$ a complex-valued matrix that depends on $t_x$, and thus on $K$, and
\begin{align}
 Q^{n} := \left(q_{i+\frac12}^n, \bar q_i^n, q_i^{(1)}, q_i^{(2)}, \ldots, q_i^{(N_1)}, q_{i,1}, q_{i,2}, \ldots, q_{i,N_2} \right) \qquad N_1,N_2 \in \mathbb N_0 
\end{align}
For the high-order Active Flux via larger stencils, $N_1 = N_2 = 0$; for high-order Active Flux via further point values, $N_1 = 0$ and for high-order via further moments, $N_2 = 0$.

$L_2$, or von Neumann stability, requires the eigenvalues of $\mathcal A$ to be contained in the closed unit disk $\{ z \in \mathbb C : |z| \leq 1\}$ for all values of $K \in [-\pi,\pi]$. With the help of the following theorem, this property can be checked more easily than by actually computing the eigenvalues:

\begin{definition}
 Given a polynomial $f \in \mathbb C[z]$, $f(z) = \sum_{j=0}^n a_j z^j$, $a_j \in \mathbb C$ with $a_n \neq 0$ and $f(0) \neq 0$, define
 \begin{align}f^*(z) &:= \sum_{j=0}^n \bar a_{n-j} z^j & f_1(z) &:= \frac{f^*(0) f(z) - f(0) f^*(z)}{z} \end{align}
  $f'$, as usual, denotes the derivative of $f$ with respect to $x$, while $\bar a$ denotes the complex-conjugate of $a \in \mathbb C$.
\end{definition}

\begin{theorem}[Schur] \label{thm:schurmiller}
 A polynomial $f \in \mathbb C[z]$ of degree $\geq 1$ has all its zeros contained in the unit disc, iff either
  \begin{itemize}
   \item $|f^*(0)| > |f(0)|$ and either $f_1$ has all its zeros contained in the unit disc or $f_1$ is a non-vanishing constant, or
   \item $f_1 \equiv 0$ and either $f'$ has all its zeros contained in the unit disc or $f'$ is a non-vanishing constant.
  \end{itemize}
\end{theorem}

Note that the recursion is sure to stop, as both $f_1$ and $f'$ are of one degree lower than $f$. This theorem was proved in \cite{schur17,schur18}; its algorithmic form was given in \cite{cohn22}, and then in \cite{miller71}, where its usefulness for stability analysis was emphasized.

Note that still this property needs to be guaranteed for all $K$. However, the Theorem \ref{thm:schurmiller} is easily amenable to an implementation which will, for any polynomial, state whether its zeros are contained in the closed unit disc up to machine error. We then propose to sample the domain of $K$, as well as any other parameters on which $\mathcal A$ might depend (in particular $\Delta t$). This allows to obtain stability results which we find very reliable in practice.

For the high-order extension of Active Flux via larger stencils, Figures \ref{fig:rk3stabilityregions3upwind1}--\ref{fig:rk3stabilityregions7thupwindminus1} show the thus obtained stability regions in the plane spanned by the CFL number and the free parameter.

\end{document}